\documentstyle{article}

\font\tenmsb=msbm10 scaled\magstep1
\font\sevenmsb=msbm7 scaled\magstep1
\font\fivemsb=msbm5 scaled\magstep1
\newfam\msbfam
\textfont\msbfam=\tenmsb
\scriptfont\msbfam=\sevenmsb
\scriptscriptfont\msbfam=\fivemsb
\def\Bbb#1{{\fam\msbfam\relax#1}}

\def\binom#1#2{{#1\choose#2}}
\def\tbinom#1#2{{\textstyle {#1\choose#2}}}

\newcommand{\cc}{{\Bbb C}}       
\newcommand{\rr}{{\Bbb R}}       
\newcommand{\zz}{{\Bbb Z}}       
\newcommand{\nn}{{\Bbb N}}       

\def\l{\ell}
\def\C{{\cal C}}
\def\Hol{{\cal H}ol}
\def\E{{\cal E}}

\def\implies{\Rightarrow}

\def\supp{\mbox{\rm supp}}

\def\iint{\int\!\!\int}

\def\a{\alpha}
\def\be{\beta}
\def\de{\delta}
\def\D{\Delta}
\def\U{\Upsilon}
\def\L{\Lambda}
\def\la{\lambda}
\def\eps{\varepsilon}
\def\g{\gamma}
\def\G{\Gamma}
\def\o{\omega}

\def\z{\zeta}
\def\s{\sigma}
\def\p{\varphi}
\def\x{\xi}
\def\u{\upsilon}

\def\zbar{\overline{\zeta}}
\def\wbar{\overline{w}}
\def\ybar{\overline{y}}
\def\xbar{\overline{x}}

\def\rr{\Bbb R}
\def\nn{\Bbb N}
\def\cc{\Bbb C}
\def\zz{\Bbb Z}
\def\tor{\Bbb T}
\def\Bbar{\overline{B^n}}

\def\ee{E\times E}
\def\eeu{(E\times E)_1}
\def\eed{(E\times E)_2}
\def\eeuu{(E\times E)_1^1}
\def\eeud{(E\times E)_1^2}

\def\unmig{\frac12}
\def\rdos{\sqrt{2}}
\def\ad{\be-\frac{n-d}{p}}
\def\as{\be-\frac{n-s}{p}}
\def\aU{\be-\frac{n-\Upsilon(E)}{p}}
\def\aL{\be-\frac{n-\Lambda(E)}{p}}
\def\jk{_j^k}
\def\lkd{_l^{k-2}}
\def\lact{{\left(1-\frac1{C\la^{3n}}\right)}}

\def\dbar{\overline{\partial}}
\def\X{\Bbb X}
\def\Y{\Bbb Y}
\def\der#1{\frac{\partial}{\partial{#1}}}

\def\Bpa{B_\a^p(\mu)}
\def\Hpb{H_\beta^p(B^n)}
\def\Lpmu{L^p(d\mu)}

\def\jet{(F_\g)_{\o(\g)<\a}}
\def\T{T^{\a}}

\def\TI{T_y^{I,\l}}

\def\dmuxy{\frac{d\mu(x)\,d\mu(y)}{\mu[x,y]^2}}

\newtheorem{teorema}{Theorem}
\newtheorem{theorem}[teorema]{Theorem}

\newtheorem{lemma}[teorema]{Lemma}

\newtheorem{proposition}[teorema]{Proposition}

\newtheorem{corolari}[teorema]{Corollary}

\newtheorem{definition}[teorema]{Definition}

\newcommand{\demo}[1]{\par\smallskip \noindent {\bf {#1}:}}
\newcommand{\qed}{\ \hfill\mbox{$\clubsuit$}}
\newcommand{\mqed}{\ \hspace*{2\fill}\clubsuit}

\newcommand{\begeq}{\begin{equation}}

\newcounter{alphacount}
\alph{alphacount}
\newenvironment{alphalist}
{
\begin{list}
{ {\rm (}{\bf \alph{alphacount}}{\rm )}}
{\usecounter{alphacount} \setlength{\rightmargin}{\leftmargin}}
}
{\end{list}}

\newcommand{\Sb}[2]{_{{\scriptstyle #1}\atop{\scriptstyle #2}}}

\title{Interpolation sets for Hardy-Sobolev spaces on the 
boundary of the unit ball of $\cc^n$.}

\author{Jaume Gudayol\thanks{Partially supported by MEC grant
PB95-0956-c02-01
and CIRIT grant GRQ94-2014. }}

\date{April 23, 1998}

\begin{document}

\maketitle

\begin{abstract}
We study the interpolation sets for the Hardy-Sobolev spaces defined
on the unit ball of $\cc^n$. 
We begin by giving a natural extension to
$\cc^n$ of the condition that is known to be necessary and suffitient 
for interpolation sets lying on the boundary of the unit disk. We show
that under this condition the restriction of a function in the
Hardy-Sobolev space to the set always exists, and lies in a Besov
space. We then show that under the assumption that there is an
holomorphic distance function for the set, there is an extension
operator from these Besov spaces to 
 the Hardy-Sobolev ones.
\end{abstract}

\section{Introduction}
In this work we study the boundary interpolation sets for
Hardy-Sobolev spaces defined on the unit ball of $\cc^n$. The study of
interpolation sets for different spaces is one
of the classical subjects of S.C.V. analysis. But in the previous works
there are serious restrictions: one considers either sets contained in
varieties or sets that have dimension less than one. In this work we
study sets not having such restrictions. Even though there are other
kinds of restrictions, we believe that one can find here a (perhaps
small) step towards the general case.

The study of interpolation sets was begun by Carleson and Rudin (See
\cite{Rudin}, chapter 10, for references). They showed
(independently) that, for $n=1$, interpolation sets for the ball algebra were
precisely those of zero Lebesgue measure. Later, and also for $n=1$,
interpolation sets for $A^\infty (D)$ were described by Alexander, Taylor and
Williams in \cite{ATW}. In this case the interpolation sets are those
satisfying that for any arc $I\subset \tor$,
$$
\frac1{|I|}\int_I \log\frac1{d(e^{it},E)}dt\le C\log\frac1{|I|}+C.
$$
Interpolation sets for the spaces $A_\a(D)$ were caracterized by
Dynkin in \cite{Dynkin-80} and Bruna in \cite{Bruna-81}.
We will say that a closed set $E\subset X$ satisfies the Uniform Hole
Condition (UHC-sets, for short) 
with respect to $X$ if there
exists $0<C<1$ so that for any
    $x\in X$ and any ball $B(x,r)$, we have
\begin{equation}
\sup\{d(y,E),\,y\in B(x,r)\}\ge C r.                 \label{holes}
\end{equation}
The UHC as
related to interpolation properties was introduced by Kotochigov, 
but other equivalent definitions have been
introduced by other authors in different contexts. The definition
says that a set has holes of a fixed size when looked at at any scale.
Dynkin and Bruna proved that, for the spaces $A_\a(D)$,
$E$ is an interpolation set iff $E$ is a UHC-set.
This characterization was obtained
by Dynkin for $\a\not\in \nn$ and by Bruna for all $0<\a<+\infty$. Later
Dynkin (\cite{Dynkin-84}) proved that a set is an interpolation set for
the Hardy-Sobolev spaces iff it is a UHC-set.

\smallskip

For $n>1$ no caracterization of boundary interpolation sets is known, not
even for the ball algebra. This does not mean that there is no information
about interpolation sets. For the ball algebra, Rudin in \cite{Rudin} devotes
all of chapter 10 to these sets, that in this case are the same as peak sets
and zero sets. There some examples are given, and one can find some background
on the problem.

Also for the ball algebra, Nagel in \cite{Nagel-76} proved that any subset of
a complex-tangential manifold is an interpolation set. On the other hand,
Davie and \O ksendal (see \cite{Rudin}, section 10.5) proved that any set that
has, in a sense, dimension less than 1 is an interpolation set. Both results
point out to the fact that an interpolation set can be as large as one wants
in the complex-tangential directions, but has to be small in the other ones.

The study of zero sets and interpolation sets for spaces other than the ball
algebra has been done by several authors. In the case of $A_\a(B^n)$ and
Hardy-Sobolev spaces results concerning sets contained in varieties were
given by Bruna and Ortega in \cite{BrunaOrtega-86}, \cite{BrunaOrtega-91},
and \cite{BrunaOrtega-93}. These works have provided us with our main
inspiration. Chaumat and Chollet obtained several results for the space
$A^\infty(B^n)$ in, for example, \cite{CC-86}, and for the Gevrey classes
in \cite{ChaumatChollet-88}.

\smallskip

Our goal was to study interpolation sets for Hardy-Sobolev spaces. But in
this case, the first problem was to know, given an $f$ in a Hardy-sobolev
space, to which function space defined on the set would the restriction
belong. This question, which is in most cases trivial, in this case is not
so. However, the results in \cite{BrunaOrtega-86} showed clearly that the
space of the restrictions should be some Besov space. But even in the real
case, no general result on restrictions of functions to Besov spaces defined
on arbitrary sets is known (however, Jonsson and Wallin in
\cite{JonssonWallin} and Jonsson in \cite{Jonsson} give some partial results).
In this article we give a restriction theorem for a general set $E\subset
S$. For the restriction to exist, we impose that the Uniform Hole
Condition (\ref{holes}) holds. We show that this condition in equivalent
to other conditions that will be useful later, and in particular, that
is equivalent to the fact that the set has, in a sense, dimension less
than the dimension of $S$.

Once we have done that, we show that under some restrictions, there is
an extension operator, thus proving that the given set is
interpolating. The restriction we impose is that we assume that there
is a holomorphic function behaving like the distance to the set. We
give some examples of such functions.



\section{Definitions and statement of results}

\subsection*{The upper dimension of a set}
 Let $(X,\rho)$ be a compact pseudo-metric space, with
 diam$(X)<+\infty$ (this means that $\rho$ satisfies the triangle
 inequality with a constant). For
$x\in X$, $R>0$ and $k\ge1$, let $N(x,R,k)$ be the maximum number of points
lying in
$B(x,kR)$ separated by a distance greater or equal than $R$.

\begin{definition}
We will say that $(X,d)\in\Upsilon_\g$ if there exists
$C(\g)=C(X,d,\g)$ such that, for all $x\in X$ and all $0<R\le kR\le1$,
$$
N(x,R,k)\le C(\g)k^\g.\eqno{(\Upsilon_\g)}
$$
We define the upper dimension $\Upsilon(X)$ as
$$
\Upsilon(X)=\inf\{\g,\,(X,d)\in\Upsilon_\g\}.
$$
\end{definition}

This dimension was first introduced by Larman under the
name of uniform metric dimension.

\smallskip

We will say that a probability measure $\mu$ lies in $U_\g=U_\g(X,\rho)$
if there exists $C(\g)$ so that for all $x\in X$ and all $0<R\le kR\le1$,
$$
\mu(B(x,kR))\le Ck^\g \mu(B(x,R)).\eqno{(U_\g)}
$$
Note that, by taking $k=1/R$, $(U_\g)$ implies the weaker condition
that, for all $x\in X$ and all $0<R\le 1$,
$$
\mu(B(x,R))\ge C R^\g.\eqno{(U'_\g)}
$$
Notice that if for some $\g$, $\mu\in U_\g$, then
$\supp\mu =X$. Moreover, in this case $\mu$ is a doubling
 measure, that is, there exists $C>0$ for which:
$$
\mu(B(x,2r))\le C\mu(B(x,r)).
$$
Let ${\cal U}=\cup_\g U_\g$. It is easily seen (see \cite{V-K}) that ${\cal
U}$ is precisely the set of all doubling measures with support on $X$.

\subsection*{The lower dimension of a set}
\begin{definition}
We will say that $(X,d)\in\Lambda_\g$ if there exists
$C(\g)=C(X,d,\g)$ such that, for $x\in X$ and $0<R\le kR\le1$,
$$
N(x,R,k)\ge C(\g)k^\g.\eqno{(\Lambda_\g)}
$$
Then we define the lower dimension $\Lambda(X)$ as:
$$
\Lambda(X)=\sup\{\g,\,(X,d)\in\Lambda_\g\}.
$$
\end{definition}
This dimension was first defined by Larman
 under the name of minimal dimension.

We will say that a doubling measure $\mu$ belongs to
$L_\g=L_\g(X,\rho)$ if
there exists $C(\g)$ so that for all $x\in X$ and all $0<R\le kR\le1$,
$$
\mu(B(x,kR))\ge Ck^\g \mu(B(x,R)).\eqno{(L_\g)}
$$
As before, by taking $k=1/R$, condition $(\Lambda_\g)$ implies
$$
\mu(B(x,R))\le C R^\g.\eqno{(L'_\g)}
$$
Note that $L_0$ poses no restriction on $\mu\in {\cal U}$.

The following improvement of Volberg and Konyagin's theorem 1 in 
\cite{V-K} can be
found in \cite{Perijo}:
\begin{theorem} \label{mutheorem}
Let $(X,\rho)$ be a pseudo-metric space. 
Let $(X,\rho)\in\Upsilon_{\u}\cap\Lambda_{\la}$, for some
$0<\la\le \u<+\infty$.
Then for any $\u'>\u$ and $\la'<\la$ (or $\la'=0$ if $\L(E)=0$) there exists
$\mu\in U_{\u'}\cap L_{\la'}$.
\end{theorem}

\subsection*{The uniform hole condition}
We are now going to restrict ourselves to closed subsets $E\subset S$. On $S$
we will use the pseudo-metric given by  $d(x,y)=|1-x\ybar|$.
In fact, we will consider $d$ as a function defined on
 $\Bbar\times\Bbar$, where $d$, although it is not a metric, satisfies
the triangle inequality
 $d(x,z)\le \sqrt{2}(d(x,y)+d(y,z))$  (for this and the following,
see \cite{Rudin}, chapter 5).

If we consider
 $\rho(x,y)=d(x,y)^\unmig$,
then $\rho$ is a metric on $S$. Then using proposition 5.1.4 in
\cite{Rudin}, we get that
 $U(S,d)=L(S,d)=n$,
and that for any $E\subset S$, $U(E,d)\le n$.

Let $\sigma$ be the normalized Lebesgue measure on $S$.
\begin{definition}
We will say that a closed set $E\subset S$ satisfies 
$\Sigma_s$ if there exists a
$C(s)$ so that for any $x\in S$, $R>0$ and $0<\eps<R$,
$$
\s(B(x,R)\cap E_\eps)\le C R^s \eps^{n-s},    \eqno{(\Sigma_s)}
$$
where $E_\eps=\{z\in S,\, d(z,E)<\eps\}$.
\end{definition}

For $n=1$, Bruna (\cite{Bruna-81}), using also results from Dynkin
(\cite{Dynkin-80}) proved that $E$ satisfied the UHC iff there was a 
$s<1$ so that
$E$ satisfied $\Sigma_s$, and that both conditions were equivalent to
 the boundedness of certain integrals. We are going to extend Dynkin
and Bruna's
results to $\cc^n$, $n \ge 1$. Namely, we are going to prove the following:

\begin{theorem} \label{longth}
Let $E\subset S$ be a closed set. Then the following statements are equivalent:
\begin{alphalist}
\item $E$ satisfies the Uniform Hole Condition;
\item There is a $s<n$ so that $(E,d)$ satisfies $\Sigma_s$;
\item $\Upsilon(E)<n$;
\item There are an $a>0$ and $C>0$ so that for any $x\in S$ and $R>0$,
$$
\int_{B(x,R)}d(z,E)^{-a}d\s(z)\le C R^{n-a};
$$
\item There is a $C>0$ so that for any $x\in S$ and $R>0$,
$$
\int_{B(x,R)}\log\left(d(z,E)^{-1}\right) d\s(z)\le
\s(B(x,R))\log\frac1R+C R^n ;
$$
\item There are $a<n$ and $C>0$ so that, for any $x\in E$ and $R>0$,
$$
\int_0^R N(x,\de,\frac{R}{\de})\de^{n-a}\frac{d\de}{\de}\le C R^{n-a};
$$
\item  There are $s<n$ and $C>0$ so that, for any $x\in E$ and $R>0$,
$$
V(B(x,R)\cap E_\eps)\le C R^{s}\eps^{n+1-s};
$$
\item There is an $a_0>1$ so that for any $a<a_0$ there is a $C>0$
 so that for any $x\in S$ and $R>0$,
$$
\int_{B(x,R)}d(z,E)^{-a}dV(z)\le C R^{n+1-a}.
$$
\end{alphalist}
\end{theorem}

\subsection*{The restriction theorem}
Let $E\subset S=\partial B^n$ be a closed set, and assume
$\Upsilon(E)<n$.
We know by theorem \ref{mutheorem} that for some
$0\le d\le s<n$ there is a measure $\mu\in U_s\cap L_d$.
Note that if $U_d\cap L_d\ne \emptyset$ for some $d$, 
then $E$ is Alhfors regular.

We will work with functions $f\in\Hpb$. Recall that if
$f\in \Hol(B^n)$, we can define its radial derivative as
$N f(z)=\sum z_j\der{z_j}f(z)$. If $f=\sum_k f_k$ is the
homogeneous expansion of $f$, we consider the fractional derivative
$$
R^\be f =\sum_k (k+1)^\be f_k
$$
so that $R^1=I+N$. Then for $p\ge 1$ and $\be>0$ the Hardy-Sobolev space
$\Hpb$ consists of those holomorphic functions such that
$$
\|f\|^p_{\Hpb}=\sup_{0<r<1}\int_S |R^\be f(r\z)|^p d\s(\z)<+\infty,
$$
where $\s$ is the normalized Lebesgue measure on $S$.

Let $X=\sum_j a_j(z)\der{z_j}+\overline{a_j(z)}\der{\overline{z}_j}$,
where $a_j\in{\cal C}^\infty(\Bbar)$, be a
vector field. We define its weight $\o(X)$ as $1/2$ if $X$ is
complex-tangential, i.e. $\sum a_j\overline{z_j} =0$, and $1$ otherwise;
 for a differential operator
$\X=X_1\cdots X_p$ define $\o(\X)=\sum \o(X_j)$. It is then known (see
\cite{AhernBruna}) that if $\o(\X)\le \be$, $\X f =X_1\cdots X_p f$ has
radial limit $\s$-a.e. on $S$, even though $\X$ may have 
order bigger than $\be$.
For such $f$ and $\X$, it makes sense to define the Hardy-Littlewood maximal
function of $\X f$ as:
$$
M(\X f)(z)=\sup_{\de>0}\frac1{V(B(z,\de)\cap\Bbar)}
\int_{B(z,\de)\cap\Bbar}|\X f(\z)|\, dV(\z),
$$
where $z\in S$.

In this context,
we want to know under which conditions there is a reasonable way of defining
the restriction  $\X f_{|E}$ and in which space of functions it lies.
For doing so, we use the results of \cite{BrunaOrtega-93}, where the following
is proved:
\begin{lemma} \label{uhuhu}
Let $\mu$ be a measure on $S$ satisfying $L_d'$, for some $d<n+1$.
Then for any differential operator  $\X$ with $\o(\X)<\ad$,
$$\int_S M\left(\X f\right)^p (\z)\,d\mu(\z)\le C \|f\|_{\Hpb}.$$
In particular, there exist $\mu$-almost everywhere the limits:
$$\lim_{r\to1}\X f(r\z),\qquad
\lim_{\de\to0}\frac1{V(B(\z,\de))}\int_{B(\z,\de)}\X f(\z)\,dV(\z),$$
and they are equal.
\end{lemma}

For a function  $f\in\Hpb$ we can define, as in
\cite{BrunaOrtega-86} and \cite{BrunaOrtega-93}, the non-isotropic Taylor
polynomial at a point $\z \in E$. This Taylor polynomial
$T^\a_\z f$ is twice as long in the complex-tangential directions.
The non-isotropic Taylor polynomial  can be defined in an
intrinsic way, using the covariant differentials of $f$, as in
\cite{BrunaOrtega-86}, or in an explicit way, using local coordinates,
as in \cite{BrunaOrtega-93}.

Let us express $\T_\z f(z)$ in coordinates:
for a point $\z\in S$, let
$w_n(z,\z)=1-z\zbar$ be the normal coordinate, and let
$w_1,\dots,w_{n-1}$ coordinate $T_\z^{\cc}$.
Then because of lemma \ref{uhuhu}, if $\g=(\g_1,\dots,\g_n)$ is a multiindex
with weight $\o(\g)=\g_n+\frac12 (\g_1+\dots+\g_{n-1})<\a$,
$$
D^\g f(\z):=\lim_{r\to1}\frac{\partial^{|\g|}}{\partial w^\g}f(r\z)=
\lim_{\de\to 0}\frac1{V(B(\z,\de))}\int_{B(\z,\de)}
\frac{\partial^{|\g|}}{\partial w^\g}f(w)\,dV(w)
$$
exists $\mu$-a.e. on $E$. Clearly $\{D^\g f(\z),\,\o(\g)<\a\}$ determines
$df(\X)(\z)$ if $\o(\X)<\a$ and we have
$$
\T_\z f(z)=\sum_{\o(\g)<\a}\frac1{\g !}D^\g f(\z)w(z,\z)^\g,
$$
where $z\in B$.

In view of that, we define the holomorphic jets of class $\Bpa$ as those
collections $F=(F_\g)_{\o(\g)<\a}$ of $L^p(\mu)$ functions such that
\begin{equation}
\sum_{\o(\g)<\a}\|F_\g\|_{\Lpmu}^p+
\sum_{\o(\g)<\a}\iint_{\ee} \frac{|F_\g(x)-D^\g(T_y^\a F)(x)|^p}
{d(x,y)^{(\a-\o(\g))p-d}}\dmuxy \label{defnormB}
\end{equation}
is finite, where $\mu[x,y]=\mu(B(x,d(x,y)))$.

Note that if $E$ is Alhfors regular the corresponding Besov space is given by
the norm:
$$
\|F\|_{\Bpa}=\sum_{\o(\g)<\a}\|F_\g\|_{\Lpmu}+\iint_{\ee}
\frac{|F_\g(x)-D^\g(T_y^\a F)(x)|^p}{d(x,y)^{(\a-\o(\g))p+d}}
d\mu(x)\,d\mu(y).$$

\demo{Remarks}
These Besov spaces with respect to $\mu$ where first introduced by
Dynkin in \cite{Dynkin-84}, when studying the interpolation problem for
$H^1_1$ in $\cc$,
and later extended to subsets of $\rr^n$ by Jonsson in \cite{Jonsson},
using only first differences of functions, when studying the restriction of a
Besov space on $\rr^n$ to a closed set.

 On the other hand, we would like to remark that our spaces can be seen as
spaces with variable regularity. For example let
$E=\Gamma_1\cup \Gamma_2$ where $\Gamma_1$ is a closed transverse curve
whereas $\Gamma_2$ is a closed complex-tangential one, and they are
disjoint. Then
if we let $\mu$ be the linear Lebesgue measure on $E$, 
$\mu\in L_{1/2}\cap U_1$.
On the other hand, on $\Gamma_1$ $d(x,y)\approx |x-y|$ whereas on $\Gamma_2$
$d(x,y)\approx |x-y|^{\frac12}$. Hence, on $\Gamma_1$
$$
d(x,y)^{\a p-d}\mu(B(x,d(x,y)))^2\approx |x-y|^{\a p-\unmig+2}
$$
whereas on $\Gamma_2$
$$
d(x,y)^{\a p-d}\mu(B(x,d(x,y)))^2\approx |x-y|^{\frac12(\a p-\unmig+2)},
$$
so that these spaces are of the Besov kind, but they have different
regularity in $\G_1$ and $\G_2$.

Another useful fact, that we will use later without further comment,
is that $\mu[x,y]\approx\mu[y,x]$.
This is so because
$$
\mu(B(x,d(x,y)))\le \mu(B(y,4 d(x,y)))\le C  \mu(B(y,d(x,y)))
$$
because of $U_s$, and if we exchange $x$ for $y$ we get the reverse
inequality. Hence our definition is symetric with respect to $x$ and $y$.

\smallskip

We have seen that for $f\in\Hpb$ and $\a=\ad$,
there is a natural way of defining the restriction
$D^\g f_{|E}$ whenever $\o(\g)<\a$.
In this case, we could ask ourselves whether $f\in\Bpa$.
The answer is yes at least if $E$ is Ahlfors-regular, with $2\a\notin\nn$,
or if $\U(E)$ and $\L(E)$ are close enough (depending on $p$).
More precisely, the result is as follows:

\begin{theorem} \label{teorestric}
Let  $E\subset S$ be a closed set with $\U(E)<n$. 
Assume that between $\aL$ and $\aU$ lies
no integer multiple of $\unmig$.
Take $n>s \ge \U(E)$ and $d \le \L(E)$ close enough 
so that the same is true for
$\ad$ and $\as$. Take any $\mu\in U_s\cap L_d$. Then for any
$f\in\Hpb$ its restriction to $E$ lies in  $\Bpa$, where
$\a=\ad$.
\end{theorem}

\noindent {\bf Remarks:}
If $E$ is a transverse curve, then the Lebesgue measure on $E$ lies in
$U_1\cap L_1$, and for these curves we recover the results in
\cite{BrunaOrtega-91}. If $E$ is a complex-tangential submanifold of real
dimension $d$, then the Lebesgue measure on this submanifold lies in
$U_{d/2}\cap L_{d/2}$, so that we recover the results in \cite{BrunaOrtega-93}.

On the other hand, for a general curve $\Gamma$ we get better results than
in \cite{BrunaOrtega-91}, because there in that case one gets the same space
as the one for a transverse curve, whereas we get a restriction theorem into
a space of variable regularity. Namely, they get that the restriction is
in a space defined by \ref{defnormB} but with the metric $|x-y|$, 
while our spaces
are defined by $d(x,y)\ge C |x-y|$. In particular, our spaces are included
into the isotropic ones.

\smallskip

The proof of the theorem is based on the representation of $f$ as an integral
of $R^\be f$, together with the use of the Bergman kernels, and
the development of $(1-z\zbar\,)^{-N}$ in a suitable way, plus the bounding
of certain integrals.

The restriction on $\ad$ and $\as$ is more or less natural. Notice that our
spaces are defined using only first differences. In the case of an
Ahlfors-regular set, the restriction is that $2\a\not\in\nn$, which is the
natural one in this case. Thus our restriction is related to the use of first
differences.

\subsection*{The extension theorem}

Consider $\Hpb$ and $\Bpa$, with $\a=\be-\frac{n-d}{p}$, defined
as above. We want to prove that, in some cases, for each jet
 $\jet\in\Bpa$ there exists $f\in\Hpb$ so that, for
 $\o(\g)<\a$, $(D^\g f)_{|E}=F_\g$ in $\Lpmu$.
In this chapter we will introduce a condition under which it holds:

\begin{definition} \label{deffundist}
Let $E\subset S$ be a closed set, and let
$h\in {\cal H}ol(B)\cap{\cal C}^{\infty}(\Bbar\setminus E)$. We will say that
$h$ is a holomorphic distance function for $E$ if:
\begin{enumerate}
\item there is a 
$C_1>0$  so that
$$
C_1^{-1} d(z,E)\le |h(z)|\le C_1 d(z,E)
$$
for all $z\in\Bbar\setminus E$;
\item for any differential
operator $\X$ there is a constant $C(\X)$ so that
$$
|\X h(z)|\le C d(z,E)^{1-\o(\X)}
$$
for all $z\in\Bbar\setminus E$.
\end{enumerate}
\end{definition}

We will give some examples of such functions in the following subsection.

\smallskip

To prove the following theorem we will work with the
Triebel-Lizorkin norms, instead of the Hardy-Sobolev ones. Let
$0<p<+\infty$, $0<q<+\infty$, and $\be\ge0$. Then the Triebel-Lizorkin
space $HF^{p,q}_\be(B^n)$ is the set of holomorphic functions $f$ on $B^n$
so that
$$
\|f\|_{p,q,\be}^p=\int_S\left(\int_0^1 (1-r^2)^{([\be]+1-\be)q-1}
|R^{[\be]+1}f(rz)|^q dr\right)^{\frac{p}{q}}d\s(z)<+\infty.
$$
It is well known that $HF^{p,2}_\be (B^n)=H^p_\be(B^n)$, and also that
$HF^{p,q_1}_\be \subset HF^{p,q_2}_\be$ if $q_2\ge q_1$.

There are two reasons for working with the Triebel-Lizorkin norms. The first
one is that it is technically simpler to work with integer powers of $R$ than
to work with $R^\be$, for $\be\not \in\nn$. 
 On the other hand, the results we get are more
general.

\begin{theorem} \label{teo_extensio}
Let $E\subset S$ be a closed set, with $\Upsilon(E)<n$.
Assume that between
$\beta-\frac{n-\Lambda(E)}p$
and $\beta-\frac{n-\Upsilon(E)}p$ lies no integer multiple of $\unmig$.
Let $d\le\Lambda(E)$ and $n>s\ge\Upsilon(E)$ be close enough so 
that this fact is still true for
$\beta-\frac{n-d}p$
and $\beta-\frac{n-s}p$, and take $\mu\in L_d\cap U_s$. Let
$\a=\beta-\frac{n-d}p$.
Assume that there is a holomorphic distance function $h$ for $E$.
Then for each  jet $\jet\in\Bpa$ there is an $f\in HF^{p,1}_\be$ so
that, for
$\o(\g)<\a$, we have $(D^\g f)_{|E}=F_\g$ in $\Lpmu$.
\end{theorem}

This theorem gives us directly that $E$ is an interpolation set for
the Hardy-Sobolev spaces.

To prove the theorem, we will first construct a function $g$ satisfying
the required growth and interpolation properties, and then we will
correct it using a $\dbar$ process to get the holomorphic function we are
looking for.

On the other hand, it is easily checked that the condition that between
$\a=\ad$ and $\as$ lies no integer multiple of $\unmig$
is needed only  to see that the function we construct lies in
$HF^{p,1}_\be$,
 and only to deal with the derivatives.
Therefore if $\a<1/2$ it is not needed. 

Observe that a theorem similar to theorem \ref{teo_extensio}, but
involving the spaces $ A_\a(B)=\Hol(B)\cap \C^\a(\Bbar)$, had already
been proved.
For these spaces, Bruna and Ortega in \cite{BrunaOrtega-86} gave the
following:
\begin{theorem}[Bruna-Ortega]
Let $\Gamma$ be a transverse curve, and let $E\subset\Gamma$ with
$\Upsilon(E)<1$. Then $E$ is an interpolation set for $A_\a$, for
$\a\in \rr^{+}\setminus \nn$.
\end{theorem}
In reading the proof, it is easy to check that the fact that $E$ is contained
in a transverse curve is used at two points of it: when, in theorem 4.3,
it is proved that for such a set there is a holomorphic distance function;
and in lemma 5.7, where it is proved that for such a set, condition
{\bf (d)} in theorem \ref{longth} is satisfied.
But in the proof of the theorem what is used is that $\U(E)<n$.
Hence, what is proved there is that under the same hipothesis as in
theorem \ref{teo_extensio}, $E$ is an interpolation set for $A^\a(B)$.

Consider now the Besov spaces $A^p_{q,\a}(B^n)$.  
It is well kwown that if we have two pairs $(q_1,\a_1)$ and $(q_2,\a_2)$
then $A^p_{q_1,\a_1}=A^p_{q_2,\a_2}$ whenever $\a_2-\a_1=(q_2-q_1)/p$.
Another remarkable fact about these spaces is that, in a limit sense,
$A^p_{0,\a}=H^p_\a$.

Let $M=\{z\in B^{n+1},\, z_{n+1}=0\}\approx B^n$.
 As $M=B^n$, we have on $M$ the spaces
$A^p_{q,\a}(M)$. We consider on $B^{n+1}$ the spaces $H^p_\be(B^{n+1})$.
Beatrous in \cite{Beatrous} proves that there exists a bounded restriction
operator
$R:H^p_\be(B^{n+1})\to A^p_{q,\a}(B^n)$ for
$\be=\a-(q-1)/p$, and that in this case there is also an extension operator
$E:A^p_{q,\a}(B^n)\to H^p_\be(B^{n+1})$ such that $R\circ E=Id$.
Therefore, using theorems \ref{teorestric} and \ref{teo_extensio},
 we can obtain similar results for these Besov spaces. Moreover, in
 the process of passing from $\cc^n$ to $\cc^{n+1}$ we drop the
 condition $\U(E)<n$.

\subsection*{Examples of holomorphic distance functions}

Here we are going to give some examples of sets $E$ for which
there is a holomorphic distance function.

\paragraph*{The Chaumat-Chollet example:}
Our first example is the
one given by Chaumat and Chollet in \cite{ChaumatChollet-88}, where they
construct a holomorphic distance function on $E$ whenever $\Upsilon(E)<1$.

They proceed as follows: they begin with any set satisfying condition {\bf
(f)}
in theorem \ref{longth} for some $n-1<a<n$; then they take for each $k$ a $2^{-k}$
covering of $E$ by balls $\{ B(\z_{j,k},2^{-k}),\,j=1,\dots,N_k\}$. Then, by
defining
$$
\phi (z)=\sum_{k=1}^\infty\sum_{j=1}^{N_k}\frac{2^{-k(n-a)}}
{2^{-k}+(1-z\zbar_{j,k})}
$$
they get a holomorphic function $\phi$ such that 
$|\phi(z)|\approx d(z,E)^{a-n}$
 and $\Re (\phi(z)) >0$. Hence
$h(z)=\phi(z)^{\frac1{a-n}}$ is the desired function. But
 saying that there is an $a$ with $n-1<a<n$ that
satisfies {\bf (f)} in theorem \ref{longth} is the same as saying that
$\U(E)<1$.
Thus we have the following analogue of Davie-\O ksendal
theorem:
\begin{corolari}  \label{coloDO}
Any set $E$
with dimension $\Upsilon(E)<1$ is an interpolation set for $\Hpb$.
\end{corolari}

We would like to remark that there is a simpler  way of constructing  this
distance function.
Let $\Upsilon(E)\le s<1$, and take $\mu\in U_s(E)$. Let $s<q_1<q_2<1$. Define
$$
h_q(z)=\int_E\frac1{(1-z\zbar\,)^q}d\mu(\z),
$$
with $q=q_1,q_2$. Then as $q<1$, we have that $\Re (1-z\zbar\,)^q\ge
C_q |1-z\zbar\,|^q$. Using it, is rather simple to see that
$|h_q(z)|\approx d(z,E)^{-q}\mu(B_z)$, and that if we take
$h(z)=h_{q_2}(z)/h_{q_1}(z)$, then $|h(z)|\approx
d(z,E)^{q_2-q_1}$. As $h$ takes values on a sector not containing 
the line $\{ \Re z<0,\, \Im z=0\}$,
we can take roots of it. Hence we can consider
$h^{1/(q_2-q_1)}$; and this is the function we were looking for.

\paragraph*{Nagel's example:}
In \cite{Nagel-76}, Alexander Nagel proves that any compact set
$K$ of a complex-tan\-gen\-tial manifold $M$ is an interpolation set for
the ball algebra. He does it by constructing a holomorphic function
with specified boundary behaviour, namely:
$$
h_p(z)=\int_M \frac{dx}{F(z,\p(x))^p}
$$
where, for $\frac{n}2<p<\frac{n}2+\frac14$, $\Re F>0$ and $|F(z,\p(x))|$
behaves as nicely so as to get that $|h_p (z)|\approx
d(z,M)^{\frac{n}2-p}$.
In particular, $M$ is an interpolation set for $\Hpb$.

If $K\subset M$ is a compact set, then it is an interpolation set for the ball
algebra, and also for $A_{\a}(B^n)$. For if we have a function $f$ on $K$, as
$M$ is totally real, we can extend it by any real method to the whole
manifold, and then extend it from the manifold to the ball.

On $\Bpa$, though, there was no known general result on the extension of
functions from subsets of $\rr^n$ to $\rr^n$. But in 
\cite{art_2} one can find
the necessary results, so that we will be able to extend any function first
to $M$, and then from $M$ to $\Hpb$. Therefore, any compact subset of
$M$ is an interpolation set for $\Hpb$, with the usual restrictions on the
indices.

\paragraph*{An interpolation set of Hausdorff dimension $n-\de$:}
For each $0<\de<1$ we can build an interpolation set with
Hausdorff dimension $n-\de$. To do so, we consider the variety
$$
\Gamma=\{z\in S,\, \Im(z_1)=\dots=\Im(z_n)=0\}.
$$
Then the Hausdorff dimension of $\Gamma$ is $n-1$, and $\Upsilon(\Gamma)=
\Lambda(\Gamma)=\frac{n-1}{2}$, because this variety is
complex-tangential.

Take $0<\de<1$ and let $C_\de\subset [-\frac12,\frac12]$
be the Cantor set with Hausdorff dimension $\de$. Then for each $t\in C_\de$
let $\G_t=e^{it}\G$ be the rotation of $\G$, that is $\{e^{it}z,\,z\in\G\}$.
Then
$$
d(\G_s,\G_t)=
   \inf\{|1-e^{i(t-s)}x\cdot y|,\, x,y\in S^{n-1}(\rr)\}
$$
and as $-1\le x\cdot y\le 1$, this minimum is of the order of
$|\sin(t-s)|\approx |t-s|$.

Let $f(z)=\frac12(1+z_1^2+\dots+z_n^2)$. Then clearly $f$ is a peak function
for $\G$, so that $f_t(z)=f(e^{-it}z)$ is a peak function for $\G_t$.
Now theorem 6.2 in \cite{BrunaOrtega-86} says that if $M$ is a
complex-tangential variety of dimension $n-1$ and $f\in A^{\infty}$
 is a peak function on
$M$, the function $h(z)=1-f(z)$ satisfies $|h(z)|\approx d(z,M)$.
Then it is easily checked that $h$  is a holomorphic distance
function for $M$. Hence for each $t\in C$ the function
$1-f_t(z)$ satisfies $|1-f_t(z)|\approx d(z,\G_t)$.

Let $E=\cup_{t\in C_\de} \G_t$. We want to construct a function $h$ so that
$|h(z)|\approx d(z,E)$. 
To do so, let $\mu$ be the Hausdorff measure on $C_\de$,
and let $\de<q<1$. Then we define
$$
h_q(z)=\int_{C_\de}\frac1{(1-f_t(z))^q}d\mu(t).
$$
This function satisfies:
\begin{equation}
|h_q(z)|\approx d(z,E)^{-(q-\de)}\qquad\mbox{ and }\qquad
\Re h_q (z)>0 \qquad \forall z\in B.
\label{tere} \end{equation}
Also, for any differential operator $\X$,
$$
|\X h_q(z)|\le C(\X) d(z,E)^{\de-q-\o(\X)}.
$$
Hence if we write $h(z)=h_q(z)^{\frac1{\de-q}}$ we have built a holomorphic
distance function for $E$, so $E$ is an interpolation set for $\Hpb$.

\smallskip

We only have to check \ref{tere}, as the other inequality is proved in
the same way. To begin with, we will check the upper inequality.

 If $I$ is an interval
centered at some $t\in C_\de$, then $\mu(I)\approx |I|^\de$. Now fix $z\in B$,
let $t_0\in E$ be so that $d(z,E)=d(z,\G_{t_0})$ and for each $t\in C_\de$
let $z_t\in\G_t$ be so that $d(z,z_t)=d(z,\G_t)$.
Let $B_k$ be the set defined by
$$
B_k=\{t\in C_\de,\, d(z,\G_t)\le 2^{k}d(z,\G_{t_0})\}
$$
for $k\ge 0$, and $B_{-1}=\emptyset$. Then if $s\in B_k$, the distance
from $s$ to $t_0$ is comparable to $d(\G_s,\G_{t_0})$. Thus, and
because of the triangle inequality,
$$
|s-t_0|\le C  d(z_s,z_{t_0})\le
 C (d(z,\G_s)+d(z,\G_{t_0}))\le C 2^{k}d(z,E),
$$
so that $\mu(B_k)\le C 2^{k\de}d(z,E)^{\de}$. Now if we decompose the
integral over $E$ into the integrals over the coronae
$B_{k+1}\setminus B_{k}$, and use the previous inequality, we obtain 
$$
|h_q(z)|\le
C \sum_{k=0}^\infty 2^{-qk}d(z,E)^{-q}
\mu(B_k\setminus B_{k-1})\le
 C d(z,E)^{-q+\de}
$$
as we wanted to see.

\smallskip

For the other bound in \ref{tere}, we use that
$$
d(z,\G_t)\le d(z,e^{i(t-t_0)}z_{t_0})
\le C (d(z,E)+|t-t_0|),
$$
whence
$\{t,\, |t-t_0|\le d(z,E)\}\subset \{t, \, d(z,\G_t)\le C d(z,E)\}$.
Now, we use that $\Re (1-f_t)>0$ and $q<1$, so that we can use the bound
$\Re(1-f_t)^q\ge C_q |1-f_t|^q$. Thus we can, modulo a constant, 
 enter the modulus inside the integral. Then we can bound the integral
 by the integral over a smaller set where we can compare $d(z,E)$ with
 $d(z,\G_t)$, and obtain the result.
\qed

\subsection{proof of theorem \ref{longth}}

We begin by proving that { ({\bf a})} implies { ({\bf b})}, which is the
hardest. To prove it, we will use the following lemma, due to Sawyer and
Wheeden (\cite{Sawyer-Wheeden}):

\begin{lemma}
\label{lemmaS-W}
Let $(X,d)$ be a separable quasi-metric space, that is, $d$ satisfies
the triangle inequality with constant $A_0$. Then for $\la=8 A_0^5$ and
for any $m\in \zz$, there are points $\{x_j^k,\,\, k\ge m,\,j=1,\dots,n_j\}$
and Borel sets $\{X_j^k,\,\, k\ge m,\,j=1,\dots,n_j\}$ (where $n_j\in \nn
\cup \{\infty\}$) such that
\begin{itemize}
\item[{\bf i)}] $B(x_j^k,\la^k)\subset X_j^k\subset
  B(x_j^k,\la^{k+1})$;
\item[{\bf ii)}] for any $k\ge m$, $\cup_{j=1}^{n_j} X_j^k= X$;
\item[{\bf iii)}] given $i$, $j$, $k$, and $l$, with $m\le k\le l$,
  either $X_i^k\subset X_j^l$ or $X_i^k\cap X_j^l=\emptyset$.
\end{itemize}
\end{lemma}

We will work in $(S,\rho)$. Here $\rho$ is a metric, so we can apply the
previous lemma with any $A_0\ge1$. Note also that there is a constant $C_0$,
depending only on $n$, so that for any
ball $Q(x,R)$ and any $y\in Q(x,R)$, if $r<R$, although $Q(y,r)$ might not be
contained in $Q(x,R)$, there is a ball $Q(w,C_0 r)\subset Q(x,R)\cap Q(y,r)$.
We will write $K_0$ for the constant appearing in \ref{holes}.

We take $\lambda=\max\{8,4/(K_0 C_0)\}$. Fix $x\in S$, $R>0$ and $0<\eps<R$,
and choose $m\in\zz$ so that $\la^m \ge 2\eps/(K_0 C_0) > \la^{m-1}$,
and apply lemma \ref{lemmaS-W} with $\la$ and $m-4$ to the pseudo-metric
space $(Q(x,R),\rho)$. Let $\{X_j^k,\,\, k\ge m-4,\,j=1,\dots,n_j\}$ be
the dyadic decomposition of $Q(x,R)$ given by the lemma. Let $k_0\in\zz$
be the only integer so that $\la^{k_0-1}<2R\le \la^{k_0}$. This implies that
$X_{1}^{k_0}=Q(x,R)$ whereas $n_{k_0-2}>1$.
Write $U_{k_0-2}=Q(x,R)$. Then $U_{k_0-2}$ satisfies
\begin{enumerate}
\item $\displaystyle U_{k}=\bigcup_{X^{k}_j\subset U_{k}}
       X^{k}_j$;
\item $E_\eps\cap Q(x,R)\subset U_{k}$,
\end{enumerate}
with $k=k_0$.
Fix $k\le k_0-2$ and assume we have built $U_k$ satisfying (1) and (2).
Assume also that $\la^k \ge 2\eps/ (K_0 C_0)$. We are going to see that we
can build $U_{k-2}$ also satisfying the previous properties and with mass
less than a constant times the mass of $U_k$.

Take a $j$ so that $X_j^k\subset U_k$. Then $Q(x_j^k,\la^k)\subset X_j^k$.
As $x\jk\in Q(x,R)$ there is a ball $Q(z, C_0\la^k)\subset Q(x\jk,\la^k)
\cap Q(x,R)$. Inside $Q(z,C_0\la^k)$ there must be a ball $Q(w,K_0 C_0\la^k)$
not intersecting $E$. If $K_0C_0\la^k>2\eps$, then $Q(w,\frac{K_0C_0}2 \la^k)
\cap E_\eps=\emptyset$.

On the other hand, $X\jk$ is the union of the sets $X\lkd$ contained in it,
so that there must be a $l_0$ so that $w\in X_{l_0}^{k-2}$, and then
$\rho(w,x_{l_0}^{k-2})\le \la^{k-1}$. Hence, $X_{l_0}^{k-2}\subset
Q(x_{l_0}^{k-2},\la^{k-1})\subset Q(w,2\la^{k-1})$. But as
$\la\ge 4/(K_0 C_0)$, also $2\la^{k-1}\le K_0C_0 \la^k/2$, and then
$Q(w,2\la^{k-1})\subset Q(w, K_0 C_0 \la^k/2)$. Thus $X_{l_0}^{k-2}\cap
E_\eps=\emptyset$. If we define
$$
\tilde{X}\jk=\bigcup_{X\lkd \cap E_\eps\ne \emptyset } X\lkd
$$
we are omiting at least $X_{l_0}^{k-2}$. So by writing
$$
U_{k-2}=\bigcup_{{X^{k}_j\subset U_{k}}}\tilde{X}\jk
$$
we obtain that $U_{k-2}$ satisfies the same conditions as $U_k$.

We are going to see that, when passing from $U_k$ to $U_{k-2}$ we are
taking away at least a fixed part of the mass of $U_k$. As
$X\jk\subset Q(x\jk,\la^{k+1})$, its mass is at most $C_n \la^{(k+1)n}$,
whereas the mass of $X\lkd$ is at least $c_n \la^{(k-2)n}$, as it contains
$Q(x\lkd,\la^{k-2})$. Then $\s(X_{l_0}^{k-2})$ is at least $1/(C\la^{3n})$
of $\s(X\jk)$. As this is true for any $j$,
$$
\s(U_{k-2})\le \lact \s(U_k).
$$
We begin with $U_{k_0-2}$ and we can go through the previous process
while $\la^{k_0-2j}\ge 2\eps/(C_0 K_0)$. As $m$ is the last integer
satisfying this inequality, we can keep on doing it while $k_0-2j\ge m$,
that is, while $2j\le k_0-m$. Take $j$ to be the last one
fulfilling this inequality. For such a $j$, we can bound 
$\s(Q(x,R)\cap E_\eps)$ by $\s(U_j)$. Applying the previous bounds, we
get that:
$$
 \s (U_j)\le \lact^{j-1}\s(U_{k_0-2})
\le C\lact^{\unmig(k_0-m)}R^{2n}.
$$
On the other hand, $\la^{m-1}<2\eps/(K_0 C_0)$ and $\la^{k_0}\ge 2R$, so that
$\la^{k_0-m+1}\ge 2 K_0 C_0 R/\eps$, that is,
$$
(k_0-m)\log\la \ge \log\frac{K_0 C_0}\la+ \log\frac{R}{\eps}.
$$
If we write $s=-\unmig(\log\lact)/\log\la$, then $s>0$ and, because of
the previous inequality,
$$
\lact^{\unmig(k_0-m)} R^{2n}=\exp(-s(k_0-m)\log\la )\, R^{2n}
\le C\eps^s R^{2n-s}
$$
so that $(E,\rho)\in\Sigma_{2n-s}$, therefore
 $(E,d)\in\Sigma_{n-\frac{s}2}$.\qed

\demo{Remark} 
We want to make clear that both the $s$ for which $E\in\Sigma_s$ and the
related constant depend only on $n$ and $K_0$. This means that if we have
two UHC-sets with the same constant $K_0$, not only they are both in the
same
$\Upsilon_s$, but they satisfy the inequality with the same constant.

\smallskip

\demo{Proof of (b)$\implies$(c)}
Assume $(E,d)\in\Upsilon_s$, and let $R>\eps>0$. We can easily reduce
us to the case $x\in E$.
 Let $x_1,\dots,x_N$ be a
maximal set of points   in
$Q(x,2R)\cap E$ with $\rho(x_i,x_j)\ge \eps$ whenever $i\ne j$. Using that
$(E,\rho)\in\Upsilon_{2s}$, we have that $N\le C2^{2s}(R/\eps)^{2s}$.
Furthermore, if $y\in Q(x,R)$ and
$\rho(y,E)<\eps$ and $x_y\in E$ is such that
$\rho(y,E)=\rho(y,x_y)$, then $x_y\in Q(x,2R)$, so that
 there exists $x_i$ for which $x_y\in Q(x_i,\eps)$, and
from here:
$$Q(x,2R)\cap\{y,\,\rho(y,E)<\eps\}\subset
\cup_{i=1}^{N}Q(x_i,2\eps).$$
Therefore,
$$
\s(Q(x,R)\cap\{y,\,\rho(y,E)<\eps\})\le
\sum_{i=1}^{N}\s(Q(x_i,2\eps))\le C R^{2s}\eps^{-2s}\eps^{2n}.
$$
Hence, taking square roots, we obtain {\bf (c)}.

\demo{Proof of (c)$\implies$(b)}
If $E$ satisfies $\Sigma_s$ and for some $R>0$ and $k\ge1$
we have points $x_1,\dots,x_N$ lying in $Q(x,kR)$ with $\rho(x_i,x_j)\ge R$,
then:
$$\s(Q(x,2kR)\cap\{y,\,\rho(y,E)<R\})\le C k^{2s}R^{2s}R^{2n-2s}=
C k^{2s}R^{2n}.$$
Then again, the balls $Q(x_i,R/2)$ are mutually disjoint, so:
$$\s(Q(x,2kR)\cap\{y,\,\rho(y,E)<R\})\ge\sum_{i=1}^N
\s(Q(x_i,\frac{R}2))\ge C N R^{2n},$$
whence $N\le C k^{2s}$, hence $(E,\rho)\in\Upsilon_{2s}$
and so $(E,d)\in\Upsilon_{s}$.

The implication {\bf (c)$\implies$(g)} is obtained in essentially
 the same way.

\demo{Proof of (b)$\implies$(d)}
Fix $x$ and $R$. If $d(x,E)\ge 2\sqrt2 R$, then for any $y\in B(x,E)$ we
have $d(y,E)\ge R$, and the bound is trivial. If $d(x,E)\le 2\sqrt2 R$,
then for any $y\in B(x,R)$ we have that $d(y,E)\le 6R$. Then if we
decompose the integral we have to bound into a sum of integrals on
coronae of decreasing radii, and apply the trivial bounds to each of
these integrals, we obtain the result.

\demo{Proof of { ({\bf d})} $\implies$ {({\bf a})}}
Let $S_R(x)=\sup \{ d(y,E),\, y\in B(x,R)\}$. Then for any $z\in B(x,R)$,
$d(z,E)^{-a} \ge S_R(x)^{-a}$. Using it to get an inferior bound of
the integral in {\bf (d)} gives us the result. To see that {\bf (e)}
implies {\bf (a)} we proceed in the same way.

\demo{Proof of { ({\bf b})} $\implies$ {({\bf e})}}
Assume $d(x,E)\le 2\sqrt2 R$, the other case being trivial.
Then the descomposition of the integral into integrals over coronae
with radii $2^{-j}R$, plus the obvious bounds for each of these
integrals, gives us the result.

\demo{Proof of { ({\bf c})} $\implies$ {({\bf f})}}
Assume $(E,d)\in\Upsilon_s$, for some $s<n$. Then using it in the
integral we have to bound gives us directly the result.

\demo{Proof of { ({\bf f})} $\implies$ {({\bf c})}}
Fix $x\in E$, $R>0$ and $k\ge1$. Then the fact that 
$N(x,\de,{kR}/{\de})$ is a decreasing function of $\de$ gives the result.

\demo{Proof of { ({\bf g})} $\implies$ {({\bf h})}}
 Just like in { ({\bf b})} $\implies$ {({\bf d})},
except that in this case $a<n+1-s$, so that $a_0=n+1-s>1$.

\demo{Proof of { ({\bf h})} $\implies$ {({\bf c})}}
Fix $x\in E$, $R>0$ and $0<\eps<R$. Then if $\{t_j,\,j=1,\dots,N\}$ is a
set of points in $B(x,R)\cap E$ with $d(t_i,t_j)\ge\eps$,
$$
C R^{n+1-a}\ge 
 \sum_{j=1}^N \int_{B(t_j,\frac\eps{2\sqrt2})}\frac1{d(z,E)^a}dV(z)\ge
C N \eps^{n+1-a}
$$
so that
$N\le C (R/\eps)^{n+1-a}$
and then $(E,d)\in \Upsilon_{n+1-a}$, where $a>1$. With this statement we
have finished the proof of the theorem.\qed

\section{Technical lemmas}
The following lemmas are going to be later.
\begin{lemma} \label{cotahmu}
Let $z\in\Bbar\setminus E$, and let $a>0$. Then
$$
\int_E\frac1{d(y,z)^{a}}
\frac{d\mu(y)}{\mu(B(z,d(y,z)))}\le C
d(z,E)^{-a}.
$$
\end{lemma}

\demo{Proof}
This is immediate if we decompose the integral in a sum of integrals
over coronae and apply the trivial bounds plus $U_s$ to each of 
these integrals.

\begin{lemma} \label{lemashitu}
Let $z\in\Bbar\setminus E$, and let $a,b,c>0$. Then
$$\iint_{\ee}\frac{d(x,y)^c}{d(x,z)^ad(y,z)^b}\dmuxy\le
Cd(z,E)^{c-a-b},$$
whenever $c>s$, $c-a-b<0$, $c-a<d$ and $c-b<d$.
\end{lemma}

\demo{Proof}
We  split  $\ee$ into the sets,
$$
\eeu=\{(x,y)\in\ee,\,d(y,z)\le d(x,z)\}
$$
and its complementary.
We will bound ony the integral over $\eeu$, as the other is bounded exactly
 in the same way, changing the roles of $x$ and $y$.
Note that, in $\eeu$, $d(x,y)\le 2\sqrt2 d(x,z)$.
 We will use it to write $\eeu$ as
$\eeu=\eeuu\cup\eeud$, where:
\begin{eqnarray*}
\eeuu&=&\{(x,y)\in\ee,\,d(x,y)\le 2\rdos d(y,z)\le 2\rdos d(x,z)\};\\
\eeud&=&\{(x,y)\in\ee,\,2\rdos d(y,z)\le d(x,y)\le 2\rdos d(x,z)\}.
\end{eqnarray*}
In $\eeuu$, $d(x,z)^{-b}\le d(y,z)^{-b}$, and also
$\eeuu\subset E(y)\times E$, where
$$E(y)=\{x\in E,\,d(y,x)\le 2\rdos d(y,z)\}.$$
Because of  $U_s$,
$$\mu(B(y,d(y,z)))=\mu(B(y,d(x,y)\frac{d(y,z)}{d(x,y)}))\le
 C\mu[y,x]\left(\frac{d(y,z)}{d(x,y)}\right)^s,$$
where $\mu[y,x]=\mu(B(y,d(y,x)))$. 
Using it, and that  $d(x,z)^{-b}\le d(y,z)^{-b}$, the integral over
 $\eeuu$ can be bounded by:
$$
\int_E\frac1{d(y,z)^{a+b-s}}\int_{E(y)}\frac{d(x,y)^{c-s}}
{\mu(B(y,d(x,y)))}d\mu(x)\frac{d\mu(y)}{\mu(B(y,d(y,z)))}.
$$
Define, for $j\ge1$,
$$E_j(y)=\{x\in E,\,2^{-j}\sqrt2 d(y,z)\le d(y,x)\le 2^{-j+1}\sqrt2
d(y,z)\}.$$
Then decomposing the integral into the integral over $E_j$ and using
the obvious bounds on each of these integrals, we can bound the
previous integral by 
$$
\int_E\sum_{j=0}^\infty
\frac{\left(2^{-j}d(y,z)\right)^{c-s}\mu(B(y,2^{-j+1}\sqrt2 d(y,z)))}
{d(y,z)^{a+b-s}\mu(B(y,2^{-j}\sqrt2 d(y,z)))}\frac{d\mu(y)}{\mu(B(y,d(y,z)))}.
$$
Now this is bounded by:
$$
\int_E\frac1{d(y,z)^{a+b-c}}
\frac{d\mu(y)}{\mu(B(z,d(y,z)))}\le C d(z,E)^{-(a+b-c)}
$$
whenever $c-s>0$, and $a+b-c>0$, as we have applied lemma \ref{cotahmu}.

In $\eeud$, $d(x,z)\le 2\rdos d(x,y)$, so
$d(x,z)\approx d(x,y)$. Also $2\rdos d(y,z)\le d(x,y)$. Then, using
$L_d$,
$$\mu(B(y,d(x,y)))=\mu(B(y,d(y,z)\frac{d(x,y)}{d(y,z)}))\ge
 C\mu(B(y,d(y,z)))\left(\frac{d(x,y)}{d(y,z)}\right)^d.$$
Using it, decomposing $E(y)$ into the sets $E_j(y)$ for $j<0$, and
 proceeding as before, we obtain the result.

\begin{lemma} \label{lemashit}
Let  $a,b,c\ge0$. If $c-a-b+n+1<0$, $c-a+n+1>0$ and $c-b+n+1>0$,
then, for any $z,w\in\Bbar$:
$$\int_{B^n}\frac{(1-|\z|)^c}{|1-\zbar z|^a|1-\zbar w|^b}dV(\z)\le
C|1-z\overline{w}|^{c-a-b+n+1}.$$
\end{lemma}

\demo{Proof}
Split $B$ into:
$$
B_1=\{\z\in B,\,|1-\zbar w|\le |1-\zbar z|\}
$$
and its complementary. Clearly it is enough to bound the integral over
$B_1$, as the other one is bounded likewise.

In  $B_1$, we can assume $b\ge n+\unmig$, by changing $a$ and $b$ for
 some $a'$ and $b'$ if necessary.
Recall that in  $B_1$, $|1-\zbar z|\ge
C(|1-\wbar z|+(1-|\z|))$. Then if we write $r=|\z|$, what we have to
 bound is
$$
\int_0^1\frac{(1-r^2)^c}
{(|1-z\wbar|+1-r)^a}\int_S\frac1{|1-r\zbar w|^b}d\s(\z)\,dr. 
$$
But proposition {\bf 1.4.10} in \cite{Rudin} says that, for $b>n$,
$$
\int_S |1-r\zbar w|^{-b}d\s(\z)\le
C(1-r^2|w|^2)^{-b+n}. 
$$
If we apply this to the last integral, and then use the change of
variables
$1-r=d(z,w)t$, we get:
$$
 C|1-\wbar z|^{c+n+1-a-b}
 \int_0^{\infty}\frac{t^{c-b+n}}{(1+t)^{a}}dt\le
C|1-\wbar z|^{c+n+1-a-b},
$$
the last integral being finite if $c-b+n>-1$ and $c+n-a-b<-1$.
\qed

\begin{proposition} \label{cotillas}
Let $a>s-n-1$, $b<n$ and $a-b+n+1>0$. Then the integral
$$
\int_B \frac{d(z,E)^a}{d(\z,z)^b}dV(z)
$$
is bounded independently of $\z\in E$.
\end{proposition}

\demo{Proof}
If $a\ge 0$ this is trivial. Assume $a<0$. Then we decompose the
integral over $B$ into the integrals over 
$B(\z,2^{-j})\setminus B(\z,2^{-j-1})$. In each of these integrals, 
$d(z,\z)\approx 2^{-j}$ and, because of part {\bf (h)} in theorem
\ref{longth}, the remaining integral can be bounded by
$2^{(-j)(n+1+a)}$ whenever $a>s-n-1$. Thus our integral is bounded by
$\sum_{j\ge 0} 2^{-i(n+1+a-b)}$, which is finite whenever $n+1+a-b>0$.
\qed

The following lemma can be found in \cite{art_2}:

\begin{lemma} \label{lemashiti}
Let   $0<b<a$. There is a constant C so that for any $z\in \cc$ with
$|z|<1$,
$$\int_0^1\left(\log\frac1t\right)^{b-1}\frac1{|1-tz|^a}dt\le C
\frac1{|1-z|^{a-b}},$$
whereas if $b>a$ this integral is bounded by a constant depending only on
$a$ and $b$.
\end{lemma}

The following lemma will allow us to compute the Taylor polynomial 
of a function written as an integral representation: 

\begin{lemma} \label{lemashitii}
For an $\a$ so that $2\a\notin\zz$, $\ell=[2\a]$,
$a\in\rr$ and $x,y\in\Bbar$, if we write
$T_y^\a=T_y^{NI,\a}$, we have:
$$\displaylines{
T_y^\a \left(\frac1{(1-\zbar z)^{a}}\right)(x)=
\sum_{k=0}^{\l}\tbinom{a+k-1}{k}
\frac{[(x-y)\zbar\,]^k}{(1-\zbar y)^{a+k}}-\hfill\cr
\hfill-\sum_{k>\l/2}^{\l}\tbinom{a+k-1}{k}\sum_{j=\l-k+1}^k
\tbinom{k}{j}\frac{\left[(x\ybar-1)\zbar y\right]^j}
{(1-\zbar y)^{a+k}}
\left[(x-y)\zbar-(x\ybar-1)\zbar y\right]^{k-j}.}$$
\end{lemma}

\demo{Proof}
In order to compute the non isotropic Taylor polynomial of weight $\a$ of
a given $F$, we begin by computing the isotropic Taylor polynomial of degree
$\ell$ of $F$. We write it in terms of $N=\sum_{i=1}^n
z_i\der{z_i}$ and $Y_j=\der{z_j}-\overline{z}_j N$ (thus $\o(N)=1$ and
$\o(Y_j)=\unmig$). Then we will keep only those terms with weight less than
$\a$, and we will be done. Recall that if  we make this development at a
fixed point $y$, we must write the polynomial in terms of $(N)_y$ and
 $(Y_j)_y$ (we are using here that the values of a tensor at a point
depend only on the values of the coefficients at that point). So we must
 write it in terms of $N_y=\sum_{i=1}^n
y_i\der{z_i}_{|z=y}$ and $(Y_j)_y=\der{z_j}_{|z=y}-\ybar_j N_y$.

Let $x,y\in  S$, and $F\in{\cal H}ol(B)$. 
Let $A=\sum_{i=1}^n(x_{i}-y_{i})(Y_{i})_y$ and $B=(x-y)\ybar
N_y$. Thus  $A$ is complex  tangential whereas  $B$ is not. 
Then a straightforward
computation shows that:
$$
(\TI F)(x)=\sum_{k=0}^\l\frac1{k!} (A+B)^k F.
$$
But  $A$ and $B$ commute, as the coefficients are frozen at $y$ and partial
derivatives commute. Hence
$$
(\TI F)(x)=\sum_{k=0}^\l\frac1{k!}\sum_{j=0}^k
\binom{k}{j}A^{k-j}B^j F.
$$

We compute now the non isotropic Taylor polynomial of weight $\a$, for
$2\a\notin\nn$.
Let $\l=[2\a]$. We begin with the isotropic Taylor polynomial of degree
$\ell$, and keep those terms with weight less than $\a$. Now the weight of
 $A^{k-j}B^j$ is
$\frac{k-j}2+j=\frac{k+j}2$. So the terms we want are those with $j\le k$
 and $j<2\a-k$, so that
$j\le\l-k$. For $k\le \l/2$, the smaller of the two is $k$, so we kave to keep
all the terms, whereas for $k> \l/2$, we have to take away the terms from
 $\l-k+1$ up to $k$. Therefore,
\begin{equation}
T_y^{NI,\a} F(x)=(\TI F)(x)
-\sum_{k>\l/2}^{\l}\frac1{k!}\sum_{j=\l-k+1}^k
\binom{k}{j} A^{k-j}B^j F. 
\label{andaqueno} 
\end{equation}

Let $\z\in B$ be fixed. We want to apply the previous computations to
$F_a(z)=(1-\zbar z)^{-a}$. Note that:
$$
N_yF_a=\sum_{i=1}^n y_i\der{z_i}_{|z=y}\,\frac1{(1-\zbar z)^a}=
\sum_{i=1}^n a
\left(\frac{y_i \zbar_i}{(1-\zbar z)^{a+1}}\right)_{|z=y}=ay\zbar
F_{a+1}(y),
$$
so $BF_a=a[(x\ybar-1)y\zbar\,]F_{a+1}(y)$. Moreover,
$$(Y_i)_yF_a=(\der{z_i}_{|z=y}-\ybar_iN_y)F_a=a\zbar_iF_{a+1}(y)-
a\ybar_i\,y\zbar F_{a+1}(y),$$
thus
$$
A F_a
=a\left[(x-y)\zbar-(x\ybar-1)\zbar y\right]F_{a+1}(y).
$$
That is, when we apply $B$ to $F_a$ we get a polynomial in $x$ and $y$
 times $F_{a+1}$. From here,
$$
B^jF_a=a(a+1)\cdot\dots\cdot(a+j-1)
\left[(x\ybar-1)\zbar y\right]^j F_{a+j}(y).
$$
Analogously,
$$
A^jF_a=a(a+1)\cdot\dots\cdot(a+j-1)
\left[(x-y)\zbar-(x\ybar-1)\zbar y\right]^jF_{a+j}(y).
$$
By adding up these two things, we get:
$$
A^{k-j}B^jF_a=\frac{(a+k-1)!}{(a-1)!}
[(x\ybar-1)\zbar y]^j
[(x-y)\zbar-(x\ybar-1)\zbar y]^{k-j}F_{a+k}(y).
$$
On the other hand,
$$
(A+B)F_a=
a[(x-y)\zbar\,]F_{a+1}(y),
$$
so that:
$$
(A+B)^jF_{a}=a(a+1)\cdot\dots\cdot(a+j-1)[(x-y)\zbar\,]^j
F_{a+j}(y).
$$
Using all of this in the formula \ref{andaqueno} gives us the claim.
\qed

\begin{lemma} \label{lemashitiii}
Let $ f \in \Hpb \cap {\cal C} (\Bbar)$. Then for
$\a<\be$,  $x,y\in \Bbar$, and $r\ge0$,
$$ T_y^\a f(x)=C
\int_0^1\int_B (\log\frac1t)^{\be-1}R^\be f(\z)(1-|\z|^2)^r
T_y^\a(\frac1{(1-t\zbar z)^{n+1+r}})(x)dV(\z)dt,$$
where $C=C(n,r,\be)$.
\end{lemma}

\demo{Proof}
We know (see \cite{Ahern-88}) that, for $f\in{\cal C}(\Bbar)$ ,
for any $z\in\Bbar$,
 $$f(z)=C(\be)\int_0^1 (\log\frac1t)^{\be-1}R^\be f(tz)\,dt.$$
But theorem 7.1.4 from \cite{Rudin} says that for $g\in H^p(B)$, with $p\ge1$,
 and $r\ge0$,
$$g(z)=C(n,r)\int_B
g(\z)\frac{(1-|\z|^2)^r}{(1-z\zbar)^{n+1+r}}dV(\z).$$
Then for $\a<\be$, we can differentiate under the integral and get the
result. \qed


\begin{lemma} \label{lemashitiv}
Let $\be>0$, $\a=\be-\frac{n-d}p$, $\g<\a$, $C\ge0$, $D>0$, and
$r\ge0$. Assume $C+D>\a+(s-d)/p=\as$ and $C<\a$.
Then for $r$ large enough  ($r>(D-\a)p-n-1)$), if we write:
$$I(x,y)=\int_0^1\left(\log\frac1t\right)^{\be-1}
\int_B\frac{|f(\z)|(1-|\z|^2)^rd(x,y)^{C+D-\g}}
{d(x,t\z)^{C+n+r+1}d(y,t\z)^{D}} dV(\z)\,dt$$
we have the bound:
$$\iint_{\ee}\frac{I(x,y)^p}{d(x,y)^{(\a-\g)p-d}}\dmuxy\le
 C\|f\|_{H^p(B)}^p.$$
\end{lemma}

\demo{Proof}
Take $\de>0$ small enough  ($\de<\frac{n+r+1}{4p'}$ and $\de<\be$).
Let $\lambda=\frac{n+r+1}{p'}-2\de$.
Because of H\"older's inequality, $I(x,y)^p$ is bounded by
$$
\int_0^1(\log\frac1t)^{(\be-\de)p-1}
\int_B\frac{|f(\z)|^p(1-|\z|^2)^rd(x,y)^{(C+D-\g)p}}
{d(x,t\z)^{(C+2\de)p+n+r+1}d(y,t\z)^{(D-4\de)p}}
dV(\z)\,dt
$$
times the integral
\begin{equation}
\left(\int_0^1(\log\frac1t)^{\de p'-1}
\int_B\frac{(1-|\z|^2)^r}
{d(x,t\z)^{\lambda p'}d(y,t\z)^{4\de p'}}
dV(\z)\,dt\right)^{\frac{p}{p'}}.                          \label{gana}
\end{equation}
We want to bound the last integral using lemma \ref{lemashit}. It is
easily checked that we can apply it whenever 
$0<4\de<\frac{n+r+1}{p'}$. To bound the integral with respect to $t$
we obtain, we use that $|1-t^2a|\approx|1-ta|$ and then apply lemma
\ref{lemashiti}. In this way we see that \ref{gana} is bounded by
$d(z,E)^{-\de p}$.

We have to evaluate 
\begin{equation}
\iint_{\ee}\frac{I(x,y)^p}{d(x,y)^{(\a-\g)p-d}}\dmuxy.     \label{lapepa}
\end{equation}
If we use the bounds we have obtained, and apply Fubini's theorem,
what we have to bound is:
$$
\int_B|f(\z)|^p(1-|\z|^2)^r\int_0^1\left(\log\frac1t\right)^{(\be-\de)p-1}
J(t\z) dt\,dV(\z),
$$
where
$$
J(t\z)= \iint_{\ee}\frac{d(x,y)^{(C+D-\a-\de)p+d}}
{d(x,t\z)^{(C+2\de)p+n+r+1}d(y,t\z)^{(D-4\de)p}}\dmuxy.
$$
We want to apply lemma \ref{lemashitu} to $J(t\z)$. It is easily
checked that we can find $C$, $D$, and then $r$ large enough and
$\de$ small enough so that we can do it. In this case, \ref{lapepa} is
bounded by:
$$
\int_B|f(\z)|^p(1-|\z|^2)^r\int_0^1\left(\log\frac1t\right)^{(\be-\de)p-1}
\frac1{d(t\z,E)^{(\be-\de)p+r+1}}dt\,dV(\z).
$$
Now if $\z_0\in E$ satisfies $d(\z,E)=d(\z,\z_0)$ and $0<t<1$, also
$d(t\z,E)\approx d(t\z,\z_0)$. From this, lemma \ref{lemashiti} and
the fact that $(1-|\z|^2)\le d(\z,E)$ we bound \ref{lapepa} by
$$
 \int_B|f(\z)|^p\frac{1}{d(\z,E)}dV(\z).
$$
But, because of part {\bf (h)} of theorem \ref{longth}, $d(z,E)^{-1}dV(z)$
 is a Carleson measure, so that this last integral is bounded by 
$\| f\|_{\Hpb}$. \qed


\begin{lemma}
\label{lemmafiver}
Let $\z,\x\in E$. Write
$$
S_1=\{ z,\in B,\,d(z,\z)\le d(z,\x) \}.
$$
Let $A$, $B$, $C$, $D$, and $F$ be $\ge0$. Then the integral
\begin{equation}
\int_{S_1}d(z,E)^{-F}\int_0^1\left(\log\frac1t\right)^{D-1}
\frac{d(tz,E)^C}{d(tz,\z)^A d(tz,\x)^B}dt\,d\s(z)          \label{easter}
\end{equation}
 is bounded by
$$
K d(\z,\x)^{-(A+B-C-D+F-n)}
$$
whenever $A+B-C-D+F-n>0$, $F<n-\U(E)$, $A-C-D>0$, and $A-C-D<n-F$.
\end{lemma}

\demo{Proof}
We split $S_1$ into $S_{11}=B(\z,\frac1{2\sqrt2}d(\z,\x))$ and 
$S_{12}$ its complementary.
For $z\in S_1$, $d(z,\z)\le d(z,\x)$. Thus, and because of triangle's
inequality, 
$d(\z,\x)\le 2\sqrt2 d(tz,\x)$. 
Then, and as $d(tz,E)\le d(tz,\z)$, the part of the integral
\ref{easter} corresponding to $S_{11}$ is bounded by
$$
d(\z,\x)^{-B}\int_{S_{11}} d(z,E)^{-F}\int_0^1
\left(\log\frac1t\right)^{D-1}
\frac{1}{d(tz,\z)^{A-C}}dt\,d\s(z).
$$
Now if $A-C-D>0$ we can apply lemma \ref{lemashiti} 
 to the inner integral, and bound the last integral by
$$
 d(\z,\x)^{-B}\int_{S_{11}} d(z,E)^{-F}
\frac{1}{d(z,\z)^{A-C-D}}d\s(z).
$$
Let $B_j=B(\z,2^{-j}\frac1{2\sqrt2}d(\z,\x))$. We split this integral
into the integrals over the coronae $B_j\setminus B_{j+1}$. In each of
these integrals we bound $d(z,\z)$ by $2^{-j} d(z,\x)$ and then apply
part {\bf (d)} of theorem \ref{longth}. Thus we obtain the bound
$$
d(\z,\x)^{-(A-C-D+F-n)}
\sum_{j=0}^\infty 2^{-j(C+D-A-F+n)},
$$
and this sum is bounded whenever $C+D-A-F+n>0$, which means $A-C-D<n-F$.

For the integral over $S_1\setminus S_{11}$, we use that
$d(tz,\x)\ge d(z,\x)\ge d(z,\z)$. Hence, and again using that
$d(tz,E)\le d(tz,\z)$, the part of \ref{easter} corresponding to
$S_{12}$ is bounded by
$$
 \int_{S_{12}} \frac{d(z,E)^{-F}}{d(z,\z)^{B}}
\int_0^1\left(\log\frac1t\right)^{D-1}
\frac{1}{d(tz,\z)^{A-C}}dt\,d\s(z).
$$
As before, if $A-C-D>0$ we can apply lemma \ref{lemashiti} 
 to the inner integral, and bound the last integral by
$$
 \int_{S_{12}} \frac{d(z,E)^{-F}}{d(z,\z)^{A+B-C-D}} d\s(z).
$$
Let $B_j=B(\z,2^{j}\frac1{2\sqrt2}d(\z,\x))$. We decompose the
integral over $S_{12}$ into the integrals over 
$B_{j}\setminus B_{j-1}$. Proceeding as before, we can bound this sum
of integrals by 
$$
 d(\z,\x)^{-(A+B-C-D+F-n)}
\sum_{j=0}^\infty 2^{-j(A+B-C-D+F-n)}
$$
and this sum is bounded whenever $A+B-C-D+F-n>0$.


\section{Proof of theorem \ref{teorestric}}
Let $f\in\Hpb$. We already know that its restriction lies in  $\Lpmu$,
so we only have to see that the integrals
$$
\iint_{\ee}
\frac{|f_\g(x)-D^\g(T_y^\a f)(x)|^p}
{d(x,y)^{(\a-\o(\g))p-d}}\dmuxy
$$
are finite. But because of lemma \ref{lemashitiii}, we have that, for
$r\ge0$, and $F_a(z,\z)=(1-z\zbar)^{-a}$, the difference 
$T_y^\a f(x)-f(x)$ can be written as 
\begin{equation}
C
\int_0^1\int_B (\log\frac1t)^{\be-1}R^\be
f(\z)(1-|\z|^2)^r \left[T_y^\a (F_{N})(tx)-
F_{N}(tx) \right]dV(\z)dt, \label{tag3.1}
\end{equation}
with $N=n+1+r$.
We want to evaluate, for  $x,y\in E\subset S$, $\z\in B$, and
$F_N(\xi)=(1-\zbar \xi)^{-N}$, the difference $T_y^\a
F_N(x)-F_N(x)$. To do so, we use that the derivatives of  
$F_N (\xi)$ with respect to
 $\xi$ are the same as those of  $F_N(\xi)-F_N(x)$.
Hence:
$$
T_y^\a(F_N)(x)-F(x)=T_y^\a(F_N-F_N(x))(x).
$$
But for $N\in\nn$,
$$
F_N(\xi)-F_N(x)=
\sum_{a=1}^{N}
\frac{(\xi-x)\zbar}{(1-\zbar \xi)^{a}(1-\zbar x)^{N-a+1}}.
$$
Thus, if we expand $(\xi-x)\zbar$ as $(1-\zbar x)+(1-\zbar \xi)$,
\begin{equation}
T_y^\a (F_N-F_N(x))(x)=
\sum_{a=1}^{N}F_{N-a+1}(x)
\left[
(1-\zbar x)T_y^\a(F_a)-
T_y^\a(F_{a-1})\right](x). \label{tag3.2}
\end{equation}

\demo{Claim}
\begin{eqnarray}
\lefteqn{(1-\zbar x)T_y^{\a} \left(F_a\right)(x)-
T_y^{\a} \left(F_{a-1}\right)(x)=}\nonumber\\
&&=-\sum_{k=m}^{\l}\tbinom{a+k-1}{k}\tbinom{k+1}{\l-k}
\frac{((x\ybar-1)y\zbar)^{\l-k}( (x-y)\zbar
-(x\ybar-1)y\zbar)^{2k-\l+1}}{(1-y\zbar\,)^{a+k}}-\nonumber\\
&&-\sum_{k=m}^{\l-1}\tbinom{a+k-1}{k}
\tbinom{k}{\l-k}\frac{(
(x\ybar-1)y\zbar
)^{\l-k+1}((x-y)\zbar -(x\ybar-1)y\zbar)^{2k-\l}}{(1-y\zbar\,)^{a+k}}.
\label{tag3.3}
\end{eqnarray}
(Recall that $\l=[2\a]$.)

Assuming this claim, which we will prove later, we can proceed with the
proof of the theorem.

\smallskip

We have to bound $D^\g_x\left(T_y^{\a}F(x)-F(x)\right)$.
We use that for $a\ne0$:
$$
N_x\left((1-\zbar z)^{a}\right)
=a(x\zbar-1)\left((1-\zbar z)^{a-1}\right)_{|z=x}+a\left((1-\zbar
z)^{a-1}\right)_{|z=x},
$$
and also that:
$$
(Y_i)_x((1-\zbar z)^{a})=a(\,\zbar_i-\xbar_i)
\left((1-\zbar z)^{a-1}\right)_{|z=x}-
a\xbar_i(x\zbar-1)\left((1-\zbar z)^{a-1}\right)_{|z=x}.$$
Iterating this we get that for $\X$ as before,
\begin{equation}
|\X_x\left((1-\zbar z)^{a}\right)|\le C(\X)d(x,\z)^{a-\o(\X)}.
\label{cotax}
\end{equation}
A similar computation shows that:
\begin{equation}
|\X_x((z-y)(\zbar-\ybar)+(z\ybar-1)(1-y\zbar\,)^a)|\le
 C d(x,y)^{\frac{a}2-\o(\X)}d(y,\z)^{\frac{a}2}.        \label{cebollon}
\end{equation}

Therefore, for $f\in\Hpb$, using \ref{tag3.1}, $|\X_x(T_y^\a f-f(x))|$
can be bounded by terms like
$$
I_j(x,y)=\int_0^1\left(\log\frac1t\right)^{\be-1}
\int_B|R^\be f(\z)|(1-|\z|^2)^r |\X T_j| dV(\z)\,dt
$$
where $j=0,1$, and
$$
T_j=\frac{\left[(x\ybar-1)ty\zbar\,\right]^{\l-k+1-j}
\left[(x-y)t\zbar-(x\ybar-1)ty\zbar\,\right]^{2k-\l+j}}
{(1-tx\zbar\,)^{N-a+1}(1-ty\zbar)^{a+k}},
$$
with $1\le a\le N$, $\l=[2\a]$, $m=[\l/2]$, and $m\le k\le\l$.
We will only evaluate the term corresponding to $T_0$, as the other
one can be dealt with likewise.

If $\o=\o(\X)$, $\o_1=\o(\X')$, $\o_2=\o(\X'')$, and
$\o_3=\o(\X''')$,
\begin{eqnarray*}
\lefteqn{
\X T_0=\frac{1}{(1-ty\zbar)^{a+k}}\sum_{\o_1+\o_2+\o_3=\o}
\X'\left((1-tx\zbar\,)^{-(N-a+1)}\right)\times
}&&\\
&&\times \X''\left(\left[(x\ybar-1)ty\zbar\,\right]^{\l-k}\right)
\X'''\left(\left[(x-y)t\zbar-(x\ybar-1)ty\zbar\,\right]^{2k-\l+1}\right);
\end{eqnarray*}
and from here, using \ref{cotax} and \ref{cebollon},
$$
|\X T_0|\le
 C\sum_{\o_1+\o'=\o}
\frac{d(x,y)^{\frac{\l+1}2-\o'}}
{d(x,t\z)^{N-a+1+\o_1}d(y,t\z)^{a+\frac{\l}2-\unmig}}.
$$
Now, using that  $1\le a\le N$, these  expressions can be bounded by:
$$
C\sum_{\o_1+\o'=\o}
\frac{d(x,y)^{\frac{\l+1}2-\o'}}
{d(x,t\z)^{N+\o_1}d(y,t\z)^{\frac{\l+1}2}}+
\frac{d(x,y)^{\frac{\l+1}2-\o'}}
{d(x,t\z)^{\unmig+\o_1}d(y,t\z)^{N+\frac{\l}2}}.
$$
We want to apply lemma \ref{lemashitiv} to these expressions, for some
$r$ large enough. In the first one,  $\g=\o$, $C=\o_1\ge0$ and
 $D=\frac{\l+1}2$, so that
$C+D-\o=\frac{\l+1}2+\o_1-\o= \frac{\l+1}2-\o'$. Hence,
$C+D\ge\frac{\l+1}2>\a$, and if between $\a$ and $\as$ lies no integer
multiple of 1/2 the matching  condition is satisfied, and also
$C<\a$  trivially. In the second term, $\g=\o$,
$C=\frac{\l}2<\a$ and
$D=\unmig+\o_1$, so $C+D=\frac{\l+1}2+\o_1$ and the requirements of lemma
\ref{lemashitiv} are also fulfilled, so that:
$$
\iint_{\ee}\frac{I_0(x,y)^p}{d(x,y)^{(\a-\o)p-d}}\dmuxy\le
 C\|R^\be f\|_{H^p(B)}^p.\mqed
$$

\demo{Proof of the claim}
We know that $T_y=T_y^{I,\l}-S_a$, where, if we write $Z=(x\ybar-1)y\zbar$,
$V=(x-y)\zbar-Z$ and $T=(1-y\zbar\,)$, we have:
$$
T_y^{I,\l} \left(\frac1{(1-\zbar \xi)^{a}}\right)(x)=
\sum_{k=0}^{\l}\tbinom{a+k-1}{k}
\frac{(V+Z)^k}{T^{a+k}},
$$
and, if $m=[\l/2]$,
$$
S_a=\sum_{k=m+1}^\l\tbinom{a+k-1}{k}\sum_{j=\l-k+1}^k
\tbinom{k}{j}\frac{Z^{j}V^{k-j}}{T^{a+k}}.
$$

Now if we substract the isotropic parts we get:
$$
(T+V+Z)
\sum_{k=0}^{\l}\tbinom{a+k-1}{k}
\frac{(V+Z)^k}{T^{a+k}}  -
\sum_{k=0}^{\l}\tbinom{a+k-2}{k}
\frac{(V+Z)^k}{T^{a+k-1}}.
$$
Rearranging these terms, and using that 
$\tbinom{a}{k}-\tbinom{a-1}{k}=\tbinom{a-1}{k}$, we see that this term
is equal to
$-\tbinom{a+\l-1}{a-1}{(V+Z)^{\l+1}}{T^{-(a+\l)}}$.

On the other hand, 
$$
-(1-\zbar x)S_a+S_{a-1}=
\left(V+Z-T\right)S_a+S_{a-1}
$$
can, rearranging terms and using the properties of the binomial
coefficient, be written as
$$
\displaylines{
\sum_{k=m+1}^{\l}\tbinom{a+k-1}{k}\sum_{j=\l-k+1}^k\tbinom{k}{j}
\frac{Z^{j}V^{k-j+1}+Z^{j+1}V^{k-j}}{T^{a+k}}-\cr
-\sum_{k=m+1}^{\l}\tbinom{a+k-2}{k-1}\sum_{j=\l-k+1}^k\tbinom{k}{j}
\frac{Z^{j}V^{k-j}}{T^{a+k-1}}. }
$$

We now split from the first sum the term corresponding to  $k=\l$.
Next we change $j$ for $j-1$ in the seccond term. Thus what we have is
$$
\displaylines{
\tbinom{a+\l-1}{\l}\sum_{j=1}^\l\tbinom{\l}{j}
\frac{Z^{j}V^{\l-j+1}+Z^{j+1}V^{\l-j}}{T^{a+\l}}+\cr
+\sum_{k=m+1}^{\l-1}\tbinom{a+k-1}{k}\sum_{j=\l-k+1}^k\tbinom{k}{j}
\frac{Z^{j}V^{k-j+1}+Z^{j+1}V^{k-j}}{T^{a+k}}-\cr
-\sum_{k=m}^{\l-1}\tbinom{a+k-1}{k}\sum_{j=\l-k}^{k+1}
(\tbinom{k}{j}+\tbinom{k}{j-1})
\frac{Z^{j}V^{k+1-j}}{T^{a+k}}.
}$$
Then we split the last term and change $j$ for $j-1$ in the latter of
the two terms we obtain. We observe also that the firs term is precisely
 $(V+Z)^{\l+1}$. Thus the 
last formula can be written as
$$
\displaylines{
\tbinom{a+\l-1}{\l}\frac{(V+Z)^{\l+1}}{T^{a+\l}}
-\tbinom{a+\l-1}{\l}
\frac{V^{\l+1}+ZV^{\l}}{T^{a+\l}}+\cr
+\sum_{k=m+1}^{\l-1}\tbinom{a+k-1}{k}\sum_{j=\l-k+1}^k\tbinom{k}{j}
\frac{Z^{j}V^{k-j+1}+Z^{j+1}V^{k-j}}{T^{a+k}}-\cr
-\sum_{k=m}^{\l-1}\tbinom{a+k-1}{k}\sum_{j=\l-k}^{k}
\tbinom{k}{j}\frac{Z^{j}V^{k+1-j}}{T^{a+k}}
-\sum_{k=m}^{\l-1}\tbinom{a+k-1}{k}\sum_{j=\l-k-1}^{k}
\tbinom{k}{j}\frac{Z^{j+1}V^{k-j}}{T^{a+k}}.
}$$
Next we split in the sums the terms that are not shared. Then we
observe that some terms with are formally in the sum, are not really 
there, because it would be needed that
 $\l-m+1\le m$, so that $m\ge(\l+1)/2$.
But $m=[\l/2]$, so this cannot happen. Thanks to this, some terms
cancel each other, so that what we have is
$$
\displaylines{
\tbinom{a+\l-1}{\l}\frac{(V+Z)^{\l+1}}{T^{a+\l}}
-\tbinom{a+\l-1}{\l}
\frac{V^{\l+1}+ZV^{\l}}{T^{a+\l}}-\cr
-\sum_{k=m}^{\l-1}\tbinom{a+k-1}{k}
(\tbinom{k}{\l-k}+\tbinom{k}{\l-k-1})
\frac{Z^{\l-k}V^{2k-\l+1}}{T^{a+k}}-\cr
-\sum_{k=m}^{\l-1}\tbinom{a+k-1}{k}
\tbinom{k}{\l-k}\frac{Z^{\l-k+1}V^{2k-\l}}{T^{a+k}}
}$$
Rearranging the terms we obtain directly that:
$$
\displaylines{
\tbinom{a+\l-1}{\l}\frac{(V+Z)^{\l+1}}{T^{a+\l}}
-\sum_{k=m}^{\l}\tbinom{a+k-1}{k}\tbinom{k+1}{\l-k}
\frac{Z^{\l-k}V^{2k-\l+1}}{T^{a+k}}-\cr
-\sum_{k=m}^{\l}\tbinom{a+k-1}{k}
\tbinom{k}{\l-k}\frac{Z^{\l-k+1}V^{2k-\l}}{T^{a+k}}.
}$$
Just adding this formula with the formula corresponding to the
isotropic Taylor polynomial gives us the result.
\qed

\section{The interpolating function} \label{sectinterpfunc}
Let $q>s$. For $z\in \Bbar\setminus E$, let
\begin{equation}
h_q(z)=\int_E |1-\zbar z|^{-q}d\mu(\z).                       \label{defhq}
\end{equation}
Let $n(z,\z)=(h_q(z)(1-\zbar z))^{-1}$.
Then our real extension has the form:
\begin{equation}
\E(f)(z)=\int_E n(z,\z)\, \T_\z f(z)d\mu(\z),              \label{defext}
\end{equation}
where $\T_y$ is the non isotropic Taylor polynomial.

The following lemma about the behaviour of $h_q(z)$ can be found (for
the isotropic metric, but the proof is valid for any pseudometric) in
\cite{art_2}:
\begin{proposition} \label{cotesh}
Let $x\notin E$, and let $x_0\in E$ be such that $d(x,E)=d(x,x_0)$.
Write $B_x=B(x_0,3 d(x,E))$. Then
\begin{enumerate} 
\item[{\bf (a)}] For $a\ge 0$, $t\in E$, and $q>s+\a$, 
$$
\int_E\frac{\rho(y,t)^a}{\rho(y,x)^q}d\mu(y)
\le C \rho(t,x)^a\rho(x,E)^{-q}\mu(B_x);
$$
\item[{\bf (b)}] For any $R>0$ and any differential operator $\X$ 
there is a $C=C(\X,R)$ so that if $|x|\le R$,
$$
|\X h_q(x)|\le C \rho(x,E)^{-q-\o(j)}\mu(B_x);
$$
\item[{\bf (c)}] 
for any $R>0$ and any differential operator $\X$ 
there is a $C=C(\X,R)$ so that if $|x|\le R$,
$$
|\X \frac1{h_q(x)}|\le C \rho(x,E)^{q-\o(\X)}\mu(B_x)^{-1}.
$$
\end{enumerate}
\end{proposition}

\medskip

Next we study the behaviour of $g$. Part {\bf (a)}
of the following lemma will give us the boundedness of the $\dbar$
correction of $g$. Part {\bf (b)} says that, in a sense, $g$ has finite
Triebel-Lizorkin norm. The fact that $g$ interpolates the jet $\jet$
up to order $\a$ is checked exactly as in the proof of theorem 8 in
\cite{art_2}, so we will not repeat it.

\begin{lemma} \label{lemma4.2}
Let $\a$ and $\beta$ be as in theorem \ref{teo_extensio}.
\begin{enumerate}
\item[{\bf (a)}]
For $k>\as$ and $\o(\X)\le k$,
$$
\int_{B^n}
d(\z,E)^{(k-\be)p-1}|\X g(\z)|^p dV(\z)<+\infty,
$$
and
$$
\int_{B^n}
\frac{d(\z,E)^{(k-\be)p-\unmig}}{(1-|\z|)^{\unmig}}|\X g(\z)|^p
dV(\z)<+\infty.
$$
\item[{\bf (b)}] If $\l=[\be]+1$, then:
$$
\int_S\left(\int_0^1 (1-t^2)^{\l-\be-1}|R^\l g(tz)|dt\right)^pd\s(z)
<+\infty.
$$
\end{enumerate}
\end{lemma}

\demo{Proof}
We will only prove {\bf (a)}, as  {\bf (b)} is proved in essentially
the same way as in \cite{art_2}.

We will only prove the boundedness of the first integral.
The proof of the second inequality is essentially the same, only just
a bit more technical. The idea is that, $(1-|\z|)^a$ being integrable
for $a<1$, we can in this case obtain bounds for the integrals with
respect to $dV$ similar to those we obtain in the case $a=0$.

We will split $\X g(z)$ in the following way:
\begin{eqnarray}
\X g(z)&=&
\int_{E} n(z,\z)\, \X(T_\z^\a f(z))d\mu(z)+\nonumber\\
&+&\sum\Sb{\o(\X')+\o(\X'')=\o(\X)}{\o(\X')>0}
\int_{E} \X'(n(z,\z))\, \X''(T_\z^\a f(z))d\mu(\z). \label{cacau}
\end{eqnarray}

Let now $z\notin E$ be fixed. The term $|\X(T_\z^\a f(z))|$ can be bounded
by sums of terms like $|f_\g (\z)|$, with $\o(\g)<\a$, independently
of $z$. Then the first term
in \ref{cacau} can be bounded by sums of terms like
$$
\frac{1}{h_q(z)}
\int_{E}\frac{|f_\g(\z)|}{d(z,\z)^q} d\mu(\z).
$$
We apply H\"older's inequality to it, and then use $U_s'$ and
proposition \ref{cotesh} to estimate the integral not containing
$|f_\g(\z)|$. Then the previous integral is bounded by  
$$
\left(
\int_{E}\frac{|f_\g(\z)|}{d(z,\z)^s} d\mu(\z).
\right)^{\frac1p}.
$$
Then we raise it to the power $p$, multiply by $d(z,\z)^{p(k-\be)-1}$,
integrate over $B$ and apply Fubini's theorem. As $s<n$, and whenever
$k>\be-\frac{n-s}p$, the integral over $B$ can be bounded using
proposition \ref{cotillas}, and we are done.

\smallskip

To estimate the second term in \ref{cacau} we use that $\int
\X(n(z,\cdot))d\mu=0$ if $\o(\X)>0$. Thus this term is bounded by sums
of terms like
\begin{equation}
\int_E |\X'(n(z,\z))|\,|\X''(T^\a_\z f(z)-T^\a_\x f(z))|\,d\mu(\z), 
                                                          \label{MMT}
\end{equation}
for any (fixed) $\x\in E$. A straightforward calculation with the help
of proposition \ref{cotesh} (see \cite{art_2} for details) shows that
\begin{equation}
|\X'(n(z,\z))|\le C(\X') d(z,E)^{q-\o(\X')}d(z,\z)^{-q}\mu(B_z). 
                                                          \label{moon}
\end{equation}
Thus, if we integrate against $d(z,\x)^{-q}d\mu(\x)$, divide by
$h_q(\z)$, and apply proposition \ref{cotesh} to $h_q$, we get that
\ref{MMT} is bounded by
\begin{equation}
\frac{d(z,E)^{2q-\o(\X')}}{\mu(B_\z)^2}
\iint_{\ee} \frac{|\X''(T^\a_\z f(z)-T^\a_\x f(z))|}
{d(z,\z)^q d(z,\x)^q} d\mu(\z) d\mu(\x). 
                                                          \label{TMM}
\end{equation}
Next we use that $\T_\z f$ is a polynomial in $z$ with degree $\le2\a$, so it
has a development at  $\x$ as:
$$
T_\x^\a f(z)=\sum_{|\g|\le2\a}\frac1{\g!}D^\g(T_\x f)(z)w(z,\z)^\g.
$$
Let $\D_\g(\z,\x)=f_\g(\z)-D^\g(T_\x^\a f)(\z)$. Then
 $\T_\z f(z) -\T_\x f(z)$ can be expressed as
\begin{equation}
\sum_{\o(\g)<\a}\frac1{\g!} \D_\g(\z,\x) w(z,\z)^\g-
\sum_{\o(\g)>\a}\frac1{\g!}
D^\g(T_\x^\a f)(\z)w(z,\z)^\g.                           \label{desc}
\end{equation}
Thus \ref{TMM} can be bounded by sums of terms like:
\begin{equation}
\frac{d(z,E)^{2q-\o(\X')}}{\mu(B_\z)^2}
\iint_{\ee} \frac{|\D_\g(\z,\x)|\,|\X''(w(z,\z)^\g)|}
{d(z,\z)^q d(z,\x)^q} d\mu(\z) d\mu(\x), 
                                                          \label{rollod}
\end{equation}
with $\o(\g)<\a$, and
\begin{equation}
\frac{d(z,E)^{2q-\o(\X')}}{\mu(B_\z)^2}
\int_E \frac{|D^\g(T_\x^\a f)(\z)|\,|\X''(w(z,\z)^\g)|}
{d(z,\z)^q d(z,\x)^q} d\mu(\z) d\mu(\x), 
                                                          \label{rollocua}
\end{equation}
with $\o(\g)>\a$. 

We begin by bounding \ref{rollocua}, which is the simplest. To do so,
we simply bound $|\X''(w(z,\z)^\g)|$ by
$d(z,\z)^{\o(\g)-\o(\X'')}$, and $|D^\g(T_\x^\a f)(\z)|$ by sums of
terms of the form $|f_\g(\x)|$. Next we integrate with respect to
$\z$, applying the bounds from proposition \ref{cotesh}. From here
on, we can proceed as for the first term in \ref{cacau}. We will only
need that $\o(\g)-(\as)>0$. But we know that $\o(\g)>\a$, so that
under the hypothesis that between $\as$ and $\ad$ lies no integer
multiple of $1/2$, the previous inequality holds. 

To bound \ref{rollod}, we begin by using that $|X''(w(z,\z))^\g|$ can be
bounded by $d(z,\z)^{\o(\g)-\o(\X'')}$. Next we raise the integral to
the power $p$, and apply H\"older's inequality to it, with some
$A,B>0$ to be chosen later. The integral not containing $\D_\g(\z,\x)$
can be bounded using proposition \ref{cotesh}. Then
\ref{rollod} can be bounded by:
\begin{equation}
\frac{d(z,E)^{(2q+\o(\g)-\o(\X)-A-B)p}}{\mu(B_\z)^2}
\iint_{\ee} \frac{|\D_\g(\z,\x)|^pd\mu(\z) d\mu(\x) }
{d(z,\z)^{(q-A)p} d(z,\x)^{(q-B)p}}. 
                                             \label{cueta}
\end{equation}
We can split $\ee$ into $\eeu(z)=\{d(\z,z)\le d(\x,z)\}$ and $\eed$ its
complementary. Clearly, we only need to bound the integral over
$\eeu$, the other one being bounded likewise.
 On $\eeu$, it is easily seen that, because of $U_s$ and
$L_d$, 
$$
\mu(B_z)\ge \frac{d(z,E)^{s}}{d(z,\x)^{s-d}d(\x,\z)^d}\mu[\z,\x]
$$
(see \cite{art_2} for details).  

We use this bound in \ref{cueta}. Then we multiply \ref{cueta} by 
$d(z,E)^{(k-\be)p-1}$, integrate over $B$ and apply Fubini's theorem.
We also use that $k-\o(\X)\ge 0$. 
To finish the proof it is then enough to see that 
$$
\int_{\{z,\, (\z,\x)\in \eeu(z)\}} 
\frac{d(z,E)^{(2q-A-B-\be+\o(\g))p-1-2s}}
{d(z,\z)^{(q-A)p-2s+2d} d(z,\x)^{(q-B)p}}dV(z)
$$
is bounded by
$$
d(\z,\x)^{-((\a-\o(\g))p+d)}.
$$

To check this last assertion, we split the set over which we are
integrating into $B_1=B(\z,d(\z,\x)/4)$
and $B_2=B\setminus B_1$ (having in mind that $d(\z,z)\le d(\x,z)$). 
On $B_1$, $d(\x,z)\ge d(\x,\z)/4$, and $d(z,E)\le d(\z,\x)/4$. Then,
if we choose $A$ so that $(q-A)p-2(s-d)<n+1$, we get the result.

On $B_2$, we use that $d(\x,z)^{-1}\le d(\z,z)^{-1}$, and 
$d(z,E)\le d(z,\z)$. Again as $(\a-\o(\g))p+d>0$, we bound the
remaining integral and we have done. \qed

\section{The $\dbar$ correction}\label{dbar}

In what follows, we will follow loosely the article
\cite{BrunaOrtega-93},
where the case when $E$ is a complex-tangential variety, with $\be\in \nn$
and in the $\Hpb$ spaces,  is considered.
Our next goal is to modify $g$ by adding to it a function which is zero up to
the necessary order on $E$ so that it is possible to get a function behaving
like $g$ on $E$ but holomorphic on $B$.

\smallskip

We want to see that from  $g={\cal E}(f)$ we can get a function lying in
$HF^{p,1}_\be$ and interpolating $\jet$. To do so, we will use the kernel:
$C_N(\z,z)=\Psi_N(\z,z)C(\z,z)$, where:
$$
\Psi_N(\z,z)=\left(\frac{1-|\z|^2}{1-z\zbar}\right)^{N}
\sum_{j=0}^{n-1}c_{j,n,N}
\left[\frac{(1-|\z|^2)(1-|z|^2)}{|1-z\zbar|^2}\right]^{j}
$$
and $C(\z,z)$ is the Cauchy kernel for the ball, that is:
$$
c_n\frac{(1-z\z)^{n-1}}
{[|1-z\zbar|^2-(1-|\z|^2)(1-|z|^2)]^{n}}
\sum_{j=1}^{n}(-1)^{j-1}(\zbar_j-\overline{z}_j)
\bigwedge_{k\ne j}d\zbar_k\bigwedge_{k=1}^{n}d\z_k.$$
This kernel was introduced by P. Charpentier  in
\cite{Charpentier} to solve the $\dbar$ problem. Namely, if $\p$ is a function
with enough regularity defined on $B^n$, then the function
$$
U(z)=T_n\p (z)=\int_{B^n}C_N(\z,z)\wedge \p(\z)
$$
satisfies $\dbar U=\p$.

\smallskip

In this section we are going to prove the following lemma, that will allow
us to finish the proof of the theorem.

\begin{lemma}\label{lemmadbar}
With the same hypothesis as in theorem \ref{teo_extensio},
let $g={\cal E}(F)$, $b>0$, and $\eta=\frac{\dbar g}{h^b}$. Let
$U(z)=T_N \eta(z)$. Then for $b$
large enough, and $N$ large enough (depending on $b$),
the following is satisfied:
\begin{enumerate}
\item[{\bf (a)}] If $\l=[\be]+1$,
$$\int_S\left(\int_0^1(1-t^2)^{\l-\be-1}|R^\l
(h^bU(t\z))|\,dt\right)^p d\s(\z)<+\infty.$$
\item[{\bf (b)}] For $\mu$-almost every $\z\in E$, if $\o(\X)<\a$,
$$\lim_{\de\to0}\frac1{V(B(\z,\de))}\int_{B(\z,\de)}
\X(h^bU)(z)dV(z)=0.$$
\end{enumerate}
 \end{lemma}

With this lemma we can easily prove theorem \ref{teo_extensio},
 for if we define $f=g-h^bU$, we have
$$
\dbar f=\dbar g-h^b\dbar U=0,
$$
so  $f$ is a holomorphic function. Moreover, part {\bf (a)} of the
lemma  together with {\bf (b)} in lemma \ref{lemma4.2} imply
that  $f\in HF^{p,1}_\be$.  What is more, part {\bf (b)} says that 
$f$ behaves like $g$  near $E$,
at least up to order $\a$, and in particular interpolates $\jet$.

\subsection{Technical lemmas}

The following lemmas are going to be used to prove lemma
\ref{lemmadbar}. When $E$ is a complex-tangential variety, they are
proved in \cite{BrunaOrtega-93}. In particular, the first proposition,
which does not depend on $E$, needs no proof.

\begin{proposition} \label{prop4.5}
Let, for $a=0,1$,
$$
D_a(\z,z)=\left[\frac{1-|\z|^2}{|1-z\zbar|}\right]^N
\frac{|1-z\zbar|^{n-1}|\z-z|^a}
{(d(\z,z)^2+[(1-|\z|^2)+(1-|z|^2)]|\z-z|^2)^{n-\unmig+\frac{a}2}}.
$$
Then the following is satisfied:
\begin{enumerate}
\item If $\z$ and $z$ are near $S$,
$$
\int_{B(z,\de)}D_a(\z,z)\,dV(\z)=O(\de^{\frac12+\frac{a}2}).
$$
\item If $N-n+1\ge0$,
$$
D_a(\z,z)=O\left(\frac{(1-|\z|^2)^N}{d(\z,z)^{N+n+\frac{a}2}}\right).
$$
\end{enumerate}
\end{proposition}

With this proposition, we can prove the following:
\begin{lemma} \label{lemma4.6}
Let $i\ge0$, $1+j+i<0$, $j+N+n+1-s>0$, and $N\ge n-1$. Then for $a=0,1$:
$$
\int_{B^n} D_a(\z,z)d(\z,z)^i d(\z,E)^j dV(\z)=
O\left(d(z,E)^{i+j+1-\frac{a}2}\right).
$$
\end{lemma}

\demo{Proof}
Both inequalities are proved in the same way, by using
proposition \ref{prop4.5}, so we will
only prove the one corresponding to $D_0$.
 To do so, we will split the integral into two
parts, over $B_1=B(z,d(z,E)/2)$ and its complementary $B_2$. Over  $B_1$,
$d(\z,z)$ and $d(\z,E)$ can be bounded by $d(z,E)$, and the remaining
integral can be estimated using part (1) of proposition \ref{prop4.5},
giving us the bound $d(z,E)^{i+j+1}$.

 On the other hand, from part
(2) of  proposition \ref{prop4.5} we get that the integral over $B_2$
can be bounded by:
$$
\int_{\frac{d(z,E)}2}^\infty\frac1{t^{N+n-i+1}}
\int_{B(z,t)}d(\z,E)^{j+N}dV(\z)dt.
$$
Now if we apply part {\bf (h)} of theorem \ref{longth} to the inner
integral, the integral with respect to $t$ can be bounded by
$d(z,E)^{j+i+1}$  whenever $j+i+1<0$.
\qed

\begin{lemma}\label{lemma4.7}
If $i>0$, $1+j+i<N$, $j+n+1-s>0$, and $N\ge n-1$, then for $a=0,1$:
\begin{eqnarray*}
\int_{B^n} D_a(\z,z)d(\z,z)^id(z,E)^j dV(z)&=&
O(d(\z,E)^j(1-|\z|)^{i+1-\frac{a}2})=\\
&=&O\left(d(\z,E)^{i+j+1-\frac{a}2}\right).
\end{eqnarray*}
\end{lemma}

\demo{Proof}
This lemma is analogous to lemma \ref{lemma4.6}, just swapping $\z$ and $z$,
and with small modifications of the indices.\qed

The following proposition can be found in \cite{BrunaOrtega-93}.

\begin{proposition} \label{lemma4.8}
Let
$$
E_1(\z,z)=\left[\frac{1-|\z|^2}{|1-z\zbar|}\right]^N\frac{|1-z\zbar|^{n-1}}
{(d(\z,z)^2+[(1-|\z|^2)+(1-|z|^2)]|\z-z|^2)^{n-\frac34}}.$$
Then the following is satisfied:
\begin{enumerate}
\item If $\z$ and $z$ are near enough to  $S$,
$$
\int_{B(z,\de)}E_1(\z,rz)\,d\s(\z)=O(\de^{\unmig});
$$
\item If $N-n+1\ge0$,
$$
E_1(\z,z)=O\left(\frac{(1-|\z|^2)^N}{d(\z,z)^{N+n-\unmig}}\right).
$$
\end{enumerate}
\end{proposition}

Just as before, using this proposition we can prove the following:
\begin{lemma} \label{lemma4.9}
Let  $i\ge0$, $j\ge0$, $0>\de>s-n$, $\unmig+j+i+\de<N$, and
$N\ge n-1$. Then, for $a=0,1$:
$$
\int_S D_a(\z,rz)d(\z,rz)^id(rz,E)^j d(z,E)^\de d\s(z)=
O\left(d(\z,E)^{j+\de}
\frac{(1-|\z|)^{i+\unmig-\frac{a}2}}{|r-|\z||^{\unmig}}\right).
$$
\end{lemma}

\demo{Proof}
As in lemma \ref{lemma4.6}, we will only prove the first inequality.
Given the relation between $D_0$ and $E_1$, and as
$d(\z,rz)\ge|r-|\z||$, it is enough to see that:
$$
\int_S E_1(\z,rz)d(\z,rz)^id(rz,E)^j d(z,E)^\de d\s(z)=
O\left(d(\z,E)^{j+\de}(1-|\z|)^{i+\unmig}\right).
$$
We will split this integral in three parts. We define
$B_1=B(\z,(1-|\z|)/2)$,
$B_2=B(\z,d(z,E)/2)\setminus B_1$ and $B_3=B^n\setminus B_2$.
The integral over $B_1$ is bounded by
$$C(1-|\z|)^{i}d(\z,E)^{j+\de}\int_{B_1}E_1(\z,rz)d\s(z)\le
C d(\z,E)^{j+\de}(1-|\z|)^{i+\unmig}.$$

Because of (2) in proposition \ref{lemma4.8},
 the integral over $B_2$ can be bounded by
$$
(1-|\z|)^Nd(\z,E)^{j+\de}\int_{S\setminus B_1}
d(\z,rz)^{-N-n+\unmig+i}d\s(z)\le C(1-|\z|)^{i+\unmig}d(\z,E)^{j+\de}.
$$

In order to bound the integral over $B_3$, we use that as
$d(z,E)\le 2d(z,\z)$ also $d(rz,E)\le C d(rz,\z)$ This,
together with (2)
from proposition \ref{lemma4.8} allows us to bound the integral over
$B_3$ by:
$$
(1-|\z|)^N\int_{\frac{d(\z,E)}2}^{+\infty}
\frac1{t^{N+n-\unmig-i-j}}\int_{B(\z,t)}d(z,E)^{\de}d\s(z)dt.
$$
We bound the inner integral using part {\bf (d)} of theorem
\ref{longth}, and we have done with it.

To prove part (2) of the lemma, it is enough to use that
$$
D_1(\z,z)\le\frac1{(1-|\z|)^{\unmig}}D_0(\z,z)
$$
and use the previous result.\qed

\begin{lemma} \label{exacte}
Let $\z\in B^n$, $0<d<\unmig$. Then
$$
\int_{\unmig}^1 \frac{(1-t)^{d-1}}{|t-|\z||^{\unmig}}dt \le
C(1-|\z|)^{-\unmig+d}.
$$
\end{lemma}

\demo{Proof}
We split the integral in three parts, namely:
$$
\int_{\unmig}^1=\int^1_{|\z|+\unmig(1-|\z|)}+
\int^{|\z|+\unmig(1-|\z|)}_{|\z|-\unmig(1-|\z|)}+
\int^{|\z|-\unmig(1-|\z|)}_{\unmig}.
$$
In order to bound the first one, we use that there
$t-|\z|>\unmig(1-|\z|)$, and the remaining integral is trivially
bounded.

To bound the second integral, we use that there $1-t\ge \unmig(1-|\z|)$
and $d-1<0$, thus
$$
\int^{|\z|+\unmig(1-|\z|)}_{|\z|-\unmig(1-|\z|)}
\frac{(1-t)^{d-1}}{|t-|\z||^{\unmig}}dt \le C (1-|\z|)^{d-1}
\int^{|\z|+\unmig(1-|\z|)}_{|\z|-\unmig(1-|\z|)}
\frac{dt}{|t-|\z||^{\unmig}}.
$$
But this last integral can be explicitely calculated, whence the
result.

To bound the third integral we observe that in case we have to consider
it, $|\z|-\unmig(1-|\z|)\ge \unmig$, and this leads to $|\z|\ge\frac23$.
Then
$$
\int^{|\z|-\unmig(1-|\z|)}_{\unmig}
\frac{(1-t)^{d-1}}{|t-|\z||^{\unmig}}dt \le
\int^{|\z|-\unmig(1-|\z|)}_0
\frac{(1-|\z|+|\z|-t)^{d-1}}{(|\z|-t)^{\unmig}}dt.
$$
Again a explicit computation of this last integral gives us the 
result.\qed

\subsection{Proof of the lemma}

Part {\bf (b)} of the lemma is proved exactly as in
\cite{BrunaOrtega-93}, using the previous lemmas, so we will not repeat
it.

 To prove {\bf (a)}, let $D_0$ and $D_1$ be defined as in proposition
\ref{prop4.5}. Then as seen in \cite{BrunaOrtega-93}, $|\X U(z)|$ is
bounded by a finite sum of terms like
$$
T_a(z)=\int_{B^n}D_a(\z,z)d(\z,z)^{k-w}d(\z,E)^{-b-k+\o(\Y)}
|\Y \dbar g(\z)|\,dV(\z),
$$
for $a=0,1$,
where $\o=\o(\X)$, $\o\le k\le 2\o$ and $\o(\Y)\le k$.

Recall that what we have to bound is:
$$
\int_S\left(\int_0^1(1-t^2)^{\l-\be-1}|R^\l h^bU(tz)|dt\right)^pd\s(z).
$$
But as $h$ is a holomorphic distance function, it is enough to bound
terms like
$$
\int_S\left(\int_0^1(1-t^2)^{\l-\be-1}d(tz,E)^{b-m}
|\X U(tz)|dt\right)^pd\s(z),
$$
with $m+\o(\X)\le \l$. But as $|\X U|$ is bounded by terms like $T_0$
or $T_1$,
we have to substitute $|\X U|$ by $T_0$ or $T_1$, and bound
these integrals. We will bound the terms of the form $T_0$, the other
ones being bounded likewise. On the other hand, the integral between $0$
and $\unmig$ is trivially bounded, thus we need only to bound the
integral between $\unmig$ and 1.

Let $0<\de<\min\{\l-\be,1/2p\}$. Because of lemmas \ref{lemma4.6} and
\ref{lemashiti}, for $b$ large enough,
$$\displaylines{
\int_{\unmig}^1 (1-t^2)^{p'(\l-\be-\de)-1}
d(tz,E)^{b+\o-\o(\Y)-1-(\l-\be)p'}\times\hfill\cr
\qquad\times\int_{B^n}
D_0(\z,tz)d(\z,tz)^{k-\o}d(\z,E)^{-b-k+\o(\Y)}dV(\z)dt\le \hfill\cr
\hfill\le C\int_{\unmig}^1 (1-t)^{p'(\l-\be-\de)-1}
d(tz,E)^{-(\l-\be)p'}dt\le C d(z,E)^{-\de p'}.
}$$
Hence, because of H\"older's inequality,
$$\displaylines{
\left(\int_\unmig^1 (1-t^2)^{\l-\be-1} d(tz,E)^{b-m}
T_1(tz)dt\right)^p\le
\hfill\cr
\qquad\le C d(z,E)^{-\de p}\int_\unmig^1 (1-t)^{\de p-1}
 d(tz,E)^{b-\frac{p}{p'}(\o-\o(\Y)-1)+(\l-\be-m)p}\times
 \hfill\cr
 \hfill\times \int_{B^n}
D_0(\z,tz)d(\z,tz)^{k-\o}d(\z,E)^{-b-k+\o(\Y)}|\Y \dbar g(\z)|^p
dV(\z)dt.
}$$
If we now integrate with respect to $z\in S$ and apply Fubinni's
theorem, we get that
$$\displaylines{
\int_S \left(\int_\unmig^1 (1-t^2)^{\l-\be-1} d(tz,E)^{b-m}
T_1(tz)dt\right)^pd\s(z)\le
\hfill\cr
\qquad\le C\int_{B^n} d(\z,E)^{-b-k+\o(\Y)}|\Y \dbar g(\z)|^p
\int_\unmig^1 (1-t)^{\de p-1}dt\,dV(\z)\times\hfill\cr
\hfill\times \int_S d(z,E)^{-\de p}
  d(tz,E)^{b-\frac{p}{p'}(\o-\o(\Y)-1)+(\l-\be-m)p}
D_1(\z,tz)d(\z,tz)^{k-\o} d\s(z).
}$$
The inner integral can be bounded using lemma \ref{lemma4.9}, and the
integral with respect to $t$ using lemma \ref{exacte}. Thus
this expression is bounded by
$$\displaylines{
\int_{B^n} d(\z,E)^{p(\o(\Y)+1-\be)-1+p(\l-\o-m)}|\Y \dbar g(\z)|^p
dV(\z)\le\hfill\cr
\hfill\le \int_{B^n} d(\z,E)^{p(\o(\Y)+1-\be)-1}|\Y \dbar
g(\z)|^p dV(\z),
}$$
as $\l\ge m+\o$. Now as $\o(\Y \dbar)=\o(\Y)+1$ this integral was
bounded in lemma \ref{lemma4.2}, so we are done.


\begin{thebibliography}{MMM,MM}





\bibitem[Ahe,88]{Ahern-88} Patrick Ahern, {\em Exceptional sets for
holomorphic Sovolev functions}, Michigan Math. J. {\bf 35} (1988), 29-41.



\bibitem[A-B,88]{AhernBruna} Patrick Ahern and Joaquim Bruna, {\em Maximal
and area integral characterizations of Hardy-Sobolev spaces in the unit
ball of $\cc^n$}, Revista Matem\'atica Iberoamericana {\bf 4} (1988),
no. 1, 123-153.






\bibitem[ATW,71]{ATW} H. Alexander, B. A. Taylor and D. L. Williams, {\em The
interpolating sets for $A^\infty$}, Journal of mathematical analysis and
aplications, {\bf 36} (1971), 556-566.





\bibitem[Bea,86]{Beatrous} Frank Beatrous, {\em Estimates for derivatives
of holomorphic functions in pseudoconvex domains}, Math. Z. {\bf 191} (1986),
91-116.





\bibitem[Bru,81]{Bruna-81} Joaquim Bruna, {\em Boundary
interpolation sets for holomorphic functions smooth up to
the boundary and BMO}, Transactions of the AMS {\bf 264} (1981),
393-409.



\bibitem[B-O,86]{BrunaOrtega-86} Joaquim Bruna and Joaquim M.
Ortega, {\em Interpolation by holomorphic functions smooth up to the boundary
in the unit ball of $\cc^n$}, Math. Ann. {\bf 274} (1986), 527-575.



\bibitem[B-O,91]{BrunaOrtega-91} Joaquim Bruna and Joaquim M.
Ortega, {\em Traces on curves of Sovolev spaces of holomorphic functions},
Arkiv f\"or Matematik {\bf 29} (1991), no.1, 25-49.



\bibitem[B-O,93]{BrunaOrtega-93} Joaquim Bruna and Joaquim M.
Ortega, {\em Interpolation along manifolds in Hardy-Sobolev spaces},
preprint, 1993.



\bibitem[B-G,98]{Perijo} Per Bylund and Jaume Gudayol, {\em 
On the existence of
doubling measures with certain regularity properties}, preprint.






\bibitem[Cha,80]{Charpentier} Philippe Charpentier, {\em Formules explicites
pour les solutions minimales de l'\'equation\ $\dbar u=f$ dans la boule et
dans le polydisque de $\cc^n$}, Ann. Inst. Fourier (Grenoble) {\bf 30}
(1980), n. 4, 121-154.



\bibitem[C-C,86]{CC-86} Jacques Chaumat and Anne-Marie Chollet,
{\em Ensembles de z\'eros et d'interpolation \`a la fronti\`ere de domaines
strictement pseudoconvexes}, Ark. Mat. {\bf 24} (1986), 27-57.




\bibitem[C-C,88]{ChaumatChollet-88} Jacques Chaumat and Anne-Marie Chollet,
{\em Classes de Gevrey non isotropes et application \`a l'interpolation},
Ann. Scuola Norm. Sup. Pisa cl. Sci. (4), {\bf 15} (1988), no. 4, 615-676.





\bibitem[Dyn,80]{Dynkin-80} E. M. Dynkin, {\em Free interpolation sets for
H\"older classes}, Math. USSR Sbornik {\bf 37} (1980), 97-117.



\bibitem[Dyn,84]{Dynkin-84} E. M. Dynkin, {\em Free interpolation by functions
 with derivatives in $H^1$}, J. Soviet Math. {\bf 27} (1984), 2475-2481.










\bibitem[Gud,98]{art_2} Jaume Gudayol, {\em Extension theorems of
    Whitney type by the use of integral operators}, 
    preprint, 1998.




\bibitem[Jon,94]{Jonsson} Alf Jonsson,
{\em Besov spaces on closed subsets of $\rr^n$}, Transactions of the AMS
{\bf 41} (1994), no. 1, 355-370.



\bibitem[J-W,84]{JonssonWallin} Alf Jonsson and Hans Wallin,
{\em Function spaces on subsets of $\rr^n$}, Mathematical Reports, Vol. 2,
part 1, Harwood Academic Publishers, Chur, 1984.










\bibitem[Nag,76]{Nagel-76} Alexander Nagel, {\em Smooth zero sets and
interpolation sets for some algebras of holomorphic functions on
strictly pseudoconvex domains}, Duke Mathematical Journal  {\bf 43} (1976),
no. 2, 323-348.



\bibitem[Rud,80]{Rudin} Walter Rudin, {\em Function theory in the
unit ball of $\cc^n$}, Springer, New York, 1980.




\bibitem[S-W,92]{Sawyer-Wheeden} E. Sawyer and R. L. Wheeden, {\em
 Weighted inequalities for fractional integrals on euclidean and
homogeneous spaces}, American Journal of Mathematics {\bf 114} (1992),
813-874.




\bibitem[V-K,88]{V-K} A.L. Vol'berg and S.V. Konyagin, {\em
On measures with the doubling condition}, Math. USSR Izvestiya
{\bf 30} (1988), 629-638.







\end{thebibliography}
\end{document}